\documentclass[11pt]{article}

\usepackage{amsfonts}
\usepackage{graphicx}
\usepackage{amsmath}
\usepackage{amsthm}
\usepackage{bbm}

\usepackage[ansinew]{inputenc}

\def\R{{\mathbb R}}

\def\N{{ \mathbb{N}}}
\def\E{{\mathbb E}}
\def\P{{\mathcal P}}

\def\Q{{\mathbb Q}}

\def \F{\mathcal F}

\newtheorem{theorem}{Theorem}

\newtheorem{corollary}[theorem]{Corollary}

\newtheorem{definition}[theorem]{Definition}

\newtheorem{lemma}[theorem]{Lemma}

\newtheorem{proposition}[theorem]{Proposition}

\def\eqlaw{\overset{d}{=}}
\def \sfleche{{\mathcal S^{\downarrow}}}

\def \coag{\textit{Coag}}
\def \0{\mathbf{0}}
\def \frag{\textit{Frag}}
\def \1{\mathbf{1}}
\newcommand{\indic}[1]{\mathbbm{1}_{\{#1\}}}

\def \cc{c_{k}}
\def \nuc{\nu_{  \coag}}
\def \ce{c_{e}}
\def \nuf{\nu_{\textit{Disl}}}
\def \meros{\bold{e}}

\author{Julien Berestycki
\footnote{
    \textit{Laboratoire de Probabilit\'es et Mod\`eles Al\'eatoires}
    \textit{Universit\'e Pierre et Marie Curie et C.N.R.S. UMR 7599}
    \textit{175, rue du Chevaleret, 75013}
    \textit{Paris, France} \qquad \&
    \textit{MODAL'X UFR SEGMI, Universit\'e Paris X - Nanterre}
    \textit{200 av. de la R\'epublique, 92000 Nanterre.}
    \textbf{e-mail:}jberest@ccr.jussieu.fr
   }
   }
\begin{document}

\title{Exchangeable Fragmentation-Coalescence processes  and their equilibrium measures}
\maketitle

\begin{abstract}
We define and study a family of Markov processes with state space
the compact set of all partitions of $\N$ that we call
exchangeable fragmen\-tation-coalescence processes. They can be
viewed as a combination of exchangeable fragmentation as defined
by Bertoin and of homogenous coalescence as defined by Pitman and
Schweinsberg or Möhle and Sagitov. We show that they admit a
unique invariant probability measure and we study some properties
of their paths and of their equilibrium measure.
\end{abstract}


\textbf{A.M.S. Classification.} 60 J 25, 60 G 09.

%\vfill\eject
\section{Introduction}

Coalescence phenomena (coagulation, gelation, aggregation,...) and
their duals fragmentation phenomena (splitting, erosion, breaks
up,...), are present in a wide variety of contexts.

References as to the fields of application of coalescence and
fragmentation models (physical chemistry, astronomy, biology,
computer sciences...) may be found in Aldous \cite{aldous_survey}
-mainly for coalescence- and in the proceedings \cite{proceeding}
for fragmentation (some further references can be found in the
introduction of \cite{moi_2}). Clearly, many fragmentation or
coalescence phenomena are not ``pure" in the sense that both are
present at the same time. For instance, in the case of polymer
formation there is a regime near the critical temperature where
molecules break up and recombine simultaneously. Another example
is given by Aldous \cite{aldous_survey}, when, in his \textit{one
specific application} section, he discusses how certain liquids
(e.g., olive oil and alcohol) mix at high temperature but separate
below some critical level. When one lowers very slowly the
temperature through this threshold, droplets of one liquid begin
to form, merge and dissolve back very quickly.

\vspace{0,7cm}

It appears that coalescence-fragmentation processes are somewhat
less tractable mathematically than pure fragmentation or pure
coalescence. One of the reasons is that by combining these
processes we lose some of the nice properties they exhibit when
they stand alone, as for instance their genealogic or branching
structure. Nevertheless, it is natural to investigate such
processes, and particularly to look for their equilibrium
measures.

In this direction Diaconis, Mayer-Wolf, Zeitouni and Zerner
\cite{diaconis_mayer-wolf} considered a coagulation-fragmentation
transformation of partitions of the interval $(0,1)$ in which the
merging procedure corresponds to the multiplicative coalescent
while the splittings are driven by a quadratic fragmentation. By
relating it to the random transposition random walk on the group
of permutations, they were able to prove a conjecture of Vershik
stating that the unique invariant measure of this Markov process
is the Poisson-Dirichlet law. We would also like to mention the
work of Pitman \cite{pitman_GEM} on a closely related split and
merge transformation of partitions of $(0,1)$ as well as Durrett
and Limic \cite{durrett_limic} on another
fragmentation-coalescence process of $(0,1)$ and its equilibrium
behavior. However, a common characteristic of all these models is
that they only allow for binary splittings (a fragment that splits
creates exactly two new fragments) and pairwise coalescences.
Furthermore the rate at which a fragment splits or merges depends
on its size and on the size of the other fragments.

\vspace{0,7cm}

Here, we will focus on a rather different class of
coagulation-fragment\-ations that can be deemed
\textit{exchangeable} or \textit{homogeneous}. More precisely,
this paper deals with processes which describe the evolution of a
countable collection of masses which results from the splitting of
an initial object of unit mass. Each fragment can split into a
countable, possibly finite, collection of sub-fragments and each
collection of fragments can merge. One can have simultaneously
infinitely many clusters that merge, each of them containing
infinitely many masses.

We will require some homogeneity property in the sense that the
rate at which fragments split or clusters merge does not depend on
the fragment sizes or any other characteristic and is not time
dependent.

Loosely speaking, such processes are obtained by combining the
semi-groups of a homogenous fragmentation and of an exchangeable
coalescent. Exchangeable coalescents, or rather $\Xi$-coalescents,
were introduced independently by Schweinsberg
\cite{schweinsberg_multiple} \footnote{ Schweinsberg was extending
the work of Pitman \cite{pitman_multiple} who treated a particular
case, the so-called $\Lambda$-coalescent in which when a
coalescence occurs, the involved fragments always merge into a
single cluster.} and by Möhle and Sagitov \cite{mohle_sagitov} who
obtained them by taking the limits of scaled ancestral processes
in a population model with exchangeable family sizes. Homogeneous
fragmentations were introduced and studied by Bertoin
\cite{bertoin_homogeneous,bertoin_self-similar,bertoin_asymptot}.

The paper is organized as follows. Precise definitions and first
properties are given in Section 3. Next, we prove that there is
always a unique stationary probability measure for these processes
and we study some of their properties. Section 5 is dedicated to
the study of the paths of exchangeable fragmentation-coalescence
processes.

The formalism used here and part of the following material owe
much to a work in preparation by Bertoin based on a series of
lectures given at the IHP in 2003, \cite{cours_bertoin}.

\section{Preliminaries}

Although the most natural state space for processes such as
fragmentation or coalescence might be the space of all possible
ordered sequence of masses of fragments
$$\sfleche = \{1 \ge x_1, \ge x_2 \ge ... \ge 0, \sum_i x_i \le 1
\},$$ as in the case of pure fragmentation or pure coalescence, we
prefer to work with the space $\P$ of partitions of $\N$. An
element $\pi$ of $\P$ can be identified with an infinite
collection of blocks (where a block is just a subset of $\N$ and
can be the empty set) $\pi = (B_1,B_2,...)$ where $\cup_i B_i
=\N$, $B_i \cap B_j =\o$ when $i \neq j$ and the labelling
corresponds to the order of the least element, i.e., if $w_i$ is
the least element of $B_i$ (with the convention $\min \o =\infty$)
then $i \le j \Rightarrow w_i \le w_j.$ The reason for such a
choice is that we can discretize the processes by looking at their
restrictions to $[n]:=\{1,...,n\}.$

As usual, an element $\pi \in \P$ can be identified with an
equivalence relation by setting $$ i \overset{\pi}{\sim} j
\Leftrightarrow i \text{ and } j \text{ are in the same block of }
\pi.$$ Let $B \subseteq B' \subseteq \N$ be two subsets of $\N$,
then a partition $\pi'$ of $B'$ naturally defines a partition $\pi
= \pi'_{|B}$ on $B$ by taking $\forall i,j \in B , i
\overset{\pi}{\sim} j \Leftrightarrow i \overset{\pi'}{\sim} j,$
or otherwise said, if $\pi'=(B_1',B_2',...)$ then $\pi= (B_1'\cap
B,B_2' \cap B,...)$ and the blocks are relabelled.

Let $\P_n$ be the set of partitions of $[n].$ For an element $\pi$
of $\P$ the restriction of $\pi$ to $[n]$ is $\pi_{|[n]}$  and we
identify each $\pi \in \P$ with the sequence
$(\pi_{|[1]},\pi_{|[2]},..) \in \P_1 \times \P_2 \times ... .$ We
endow $\P$ with the distance
$$d(\pi^1,\pi^2) = 1/ \max\{ n \in \N : \pi_{|[n]}^1 =
\pi_{|[n]}^2\}.$$ The space $(\P,d)$ is then compact. In this
setting it is clear that if a family $(\Pi^{(n)})_{n \in \N}$ of
$\P_n$-valued random variable is \textit{compatible}, i.e., if for
each $n$ $$\Pi^{(n+1)}_{|[n]}=\Pi^{(n)} \; \text{a.s.},$$ then,
almost surely, the family $(\Pi^{(n)})_{n \in \N}$ uniquely
determines a $\P$-valued variable $\Pi$ such that for each $n$ one
has $$\Pi_{|[n]} = \Pi^{(n)}.$$ Thus we may define the
exchangeable fragmentation-coalescence processes by their
restrictions to $[n]$.

Let us now define deterministic notions which will play a crucial
role in the forthcoming constructions. We define two operators on
$\P$, a coagulation operator, $\pi,\pi' \in \P \mapsto
\coag(\pi,\pi')$ (the coagulation of $\pi$ by $\pi'$) and a
fragmentation operator $\pi,\pi'\in \P, k\in \N \mapsto
\frag(\pi,\pi',k)$ (the fragmentation of the $k$-th block of $\pi$
by $\pi'$).
\begin{itemize}
\item Take $\pi =(B_1,B_2,...)$ and $\pi'=(B'_1,B'_2,...)$. Then
$\coag(\pi,\pi') = (B''_1, B''_2,...),$ where $B''_1=\cup_{i \in
B'_1}B_i , B''_2=\cup_{i \in B'_2}B_i,...$  Observe that the
labelling is consistent with our convention. \item Take $\pi
=(B_1,B_2,...)$ and $\pi'=(B'_1,B'_2,...).$ Then
$\frag(\pi,\pi',k)$ is the relabelled collection of blocks formed
by all the $B_i$ for $i \neq k$, plus the sub-blocks of $B_k$
given by $\pi'_{|B_k}$.
\end{itemize}

Similarly, when $\pi \in \P_n$ and $\pi' \in \P$ or $\pi' \in
\P_k$ for $k \ge \#\pi$  (where $\#\pi$ is the number of non-empty
blocks of $\pi$) one can define $\coag(\pi,\pi')$ as above and
when $\pi' \in \P$ or $\pi' \in \P_m$ for $m \ge \text{Card}(B_k)$
one can define $\frag(\pi,\pi',k)$ as above.

Define $\0 :=(\{1\},\{2\},...)$ the partition of $\N$ into
singletons, $\0_n:=\0_{|[n]}$, and $\1:=(\{1,2,...\})$ the trivial
partition of $\N$ in a single block, $\1_n:=\1_{|[n]}$. Then $\0$
is the neutral element for $\coag,$ i.e., for each $\pi \in \P $
$$\coag(\pi,\0)=\coag(\0,\pi)=\pi,$$ (for $\pi \in \cup_{n \ge 2}
\P_n$, as $\coag(\0,\pi)$ is not defined, one only has
$\coag(\pi,\0)=\pi$) and $\1$ is the neutral element for $\frag,$
i.e., for each $\pi \in \P$ one has
$$\frag(\1,\pi,1)=\frag(\pi,\1,k)=\pi.$$ Similarly, when $\pi \in
\cup_{n \ge 2} \P_n, $ for each $k \le \# \pi$ one only has
$$\frag(\pi,\1,k)=\pi.$$

Note also that the coagulation and fragmentation operators are not
really reciprocal because $\frag$ can only split one block at a
time.

Much of the power of working in $\P$ instead of $\sfleche$ comes
from Kingman's theory of exchangeable partitions. For the time
being, let us just recall the basic definition. Define the action
of a permutation $\sigma : \N \mapsto \N$ on $\P$ by
$$ i \overset{\sigma(\pi)}{\sim} j \Leftrightarrow \sigma(i)
\overset{\pi}{\sim} \sigma(j).$$ A random element $\Pi$ of $\P$ or
a $\P$ valued process $\Pi(\cdot)$ is said to be exchangeable if
for any permutation $\sigma$ such that $\sigma(n)=n$ for all large
enough $n$ one has $\sigma(\Pi) \eqlaw \Pi$ or $\Pi(\cdot)\eqlaw
\sigma(\Pi(\cdot)).$

\section{Definition, characterization and construction of EFC
processes}

\subsection{Definition and characterization}
We can now define precisely the exchangeable
fragmentation-coalescence processes and state some of their
properties. Most of the following material is very close to the
analogous definitions and arguments for pure fragmentations (see
\cite{bertoin_homogeneous}) and coalescences (see
\cite{pitman_multiple,schweinsberg_multiple}).
\begin{definition} \label{fundamental definition}
A $\P$-valued Markov process $(\Pi(t),t \ge 0)$, is an
exchangeable fragmentation-coalescent process (``EFC process"
thereafter) if it has the following properties:
\begin{itemize}
\item It is exchangeable. \item Its restrictions $\Pi_{|[n]}$ are
càdlàg finite state Markov chains which can only evolve by
fragmentation of one block  or by coagulation.
\end{itemize}
More precisely, the transition rate of $\Pi_{|[n]}(\cdot)$ from
$\pi$ to $\pi',$ say $q_n(\pi,\pi')$, is non-zero only if $\exists
\pi''$ such that $\pi'=\coag(\pi,\pi'')$ or $\exists \pi'',k\ge 1$
such that $\pi'=\frag(\pi,\pi'',k)$.
\end{definition}

Remark that this definition implies that $\Pi(0)$ should be
exchangeable. Hence the only possible deterministic starting
points are $\1$ and $\0$ because the measures $\delta_{\1}(\cdot)$
and $\delta_{\0}(\cdot)$ (where $\delta_{\bullet}(\cdot)$ is the
Dirac mass in $\bullet$) are the only exchangeable measures of the
form $\delta_{\pi}(\cdot).$ If $\Pi(0)=\0$ we say that the process
is started from dust, and if $\Pi(0)=\1$ we say it is started from
unit mass.

Note that the condition that the restrictions $\Pi_{|[n]}$ are
càdlàg implies that $\Pi$ itself is also càdlàg.

Fix $n$ and $\pi \in \P_n$. For convenience we will also use the
following notations for the transition rates: For $\pi' \in \P_m
\backslash \{\0_m\}$ where $m = \# \pi$ the number of non-empty
blocks of $\pi$, call
$$C_n(\pi,\pi'):=q_n(\pi,\coag(\pi,\pi'))$$ the rate of
coagulation by $\pi'$. For $k \le \# \pi$ and $\pi' \in \P_{|
B_k|} \backslash \{ \1_{|B_k|} \}$ where $|B_k|$ is the cardinal
of the $k$-th block, call
$$F_n(\pi,\pi',k):=q_n(\pi,\frag(\pi,\pi',k))$$ the rate of
fragmentation of the $k$th block by $\pi'.$

We will say that an EFC process is non-degenerated if it has both
a fragmentation and coalescence component, i.e., for each $n$
there are some $\pi_1' \neq \1_n$ and $\pi'_2 \neq \0_n$ such that
$F_n(\1_n,\pi'_1,1)>0$ and $C_n(\0_n,\pi'_2)>0.$

Of course the compatibility of the $\Pi_{|[m]}$ and the
exchangeability requirement entail that not every family of
transition rates is admissible. In fact, it is enough to know how
$\Pi_{|[m]}$ leaves $\1_m$ and $\0_m$ for every $m \le n$ to know
all the rates $q_n(\pi,\pi').$

\begin{proposition} \label{prop taux de sauts}
There exists two families $((C_n(\pi))_{\pi \in \P_n \backslash
\{\0_n\}})_{n \in \N}$ and \\ $((F_n(\pi))_{\pi \in \P_n
\backslash \{\1_n\}})_{n \in \N}$ such that for every $m \le n$
 and for every $\pi \in \P_n$ with $m$ blocks ($\# \pi =m$) one
 has
\begin{enumerate}
\item   For each $\pi' \in \P_m \backslash \{\0_m\}$
$$q_n(\pi,\coag(\pi,\pi')) =C_n(\pi,\pi')= C_m(\pi').$$ \item For
each $k \le m$ and for each $\pi' \in \P_{|B_k|} \backslash
\{\1_{|B_k|}\},$
$$q_n(\pi,\frag(\pi,\pi',k))=F_n(\pi,\pi',k)=F_{|B_k|}(\pi').$$ \item All other
transition rates are zero.
\end{enumerate}

Furthermore, these rates are exchangeable, i.e., for any
permutation $\sigma$ of $[n]$, for all $\pi \in \P_n$ one has
$C_n(\pi)=C_n(\sigma(\pi))$ and $F_n(\pi)=F_n(\sigma(\pi)).$
\end{proposition}

As the proof of this result is close to the arguments used for
pure fragmentation or pure coalescence and is rather technical, we
postpone it until section 6.

\vspace{1cm}

Observe that, for $n$ fixed, the finite families $(C_n(\pi))_{\pi
\in \P_n \backslash \{\0_n\}}$ and \\ $(F_n(\pi))_{\pi \in \P_n
\backslash \{\1_n\}}$ may be seen as measures on $\P_n$. The
compatibility of the $\Pi_{|[n]}(\cdot)$ implies the same property
for the $(C_n,F_n)$, i.e., as measures, the image of $C_{n+1}$
(resp. $F_{n+1}$) by the projection $\P_{n+1} \mapsto \P_n$ is
$C_n$ (resp. $F_n$), see Lemma 1 in \cite{bertoin_homogeneous} for
a precise demonstration in the case where there is only
fragmentation ($C \equiv 0$), the general case being a simple
extension. Hence, by Kolmogorov's extension Theorem, there exists
a unique measure $C$ and a unique measure $F$ on $\P$ such that
for each $n$ and for each $\pi \in \P_n$ such that $\pi \neq \1_n$
(resp. $\pi \neq \0_n$) $$C_n(\pi)=C(\{ \pi' \in \P :
\pi'_{|[n]}=\pi \}) \text{ resp. } F_n(\pi)=F(\{ \pi' \in \P :
\pi'_{|[n]}=\pi \}).$$

Furthermore, as we have remarked,  the measures $C_n$ and $F_n$
are exchangeable. Hence, $C$ and $F$ are exchangeable measures.
They must also verify some integrability conditions because the
$\Pi_{|[n]}(\cdot)$ are Markov chains and have thus a finite jump
rate at any state. For $\pi \in \P$ define $Q(\pi,n) :=\{\pi' \in
\P : \pi'_{|[n]}=\pi_{|[n]} \}.$ Then for each $n \in \N$ we must
have
$$C(\P \backslash Q(\0,n)) <  \infty$$ and $$F(\P \backslash
Q(\1,n)) < \infty.$$ It is clear that we can suppose without loss
of generality that $C$ and $F$ assign no mass to the respective
neutral elements for $\coag$ and $\frag$, i.e., $C(\0)=0$ and
$F(\1)=0$.

Here are three simple examples of exchangeable measures.
\begin{enumerate}
\item Let $\epsilon_n$ be the partition that has only two non
empty blocks: $\N \backslash \{n\}$ and $\{n\}$. Then the
(infinite) measure $\meros(\cdot) = \sum_{n\in\N}
\delta_{\epsilon_n}(\cdot)$ (where $\delta$ is the Dirac mass) is
exchangeable. We call it the erosion measure . \item For each $i
\neq j \in \N,$ call $\epsilon_{i,j}$ be the partition that has
only one block which is not a singleton: $\{i,j\}$. Then the
(infinite) measure $\kappa(\cdot) =  \sum_{i <j \in\N}
\delta_{\epsilon_{i,j}}(\cdot)$ is exchangeable. We call it the
Kingman measure. \item Take $x \in \sfleche := \{ x_1 \ge x_2 \ge
... \ge 0 ; \sum_i x_i \le 1 \}.$ Let $(X_i)_{i \in \N}$ be a
sequence of independent variables with respective law given by
$P(X_i =k) = x_k$ for all $k \ge 1$ and $P(X_i=-i)=1-\sum_j x_j$.
Define a random variable $\pi$ with value in $\P$ by letting $ i
\overset{\pi}{\sim} j \Leftrightarrow X_i=X_j.$ Following Kingman,
we call $\pi$ the $x$-paintbox process and denote by $\mu_x$ its
distribution. Let $\nu$ be a measure on $\sfleche$, then the
mixture $\mu_{\nu}$ of paintbox processes directed by $\nu$, i.e.,
$$\mu_{\nu}(A)=\int_{\sfleche} \mu_x(A) \nu(dx),$$ is an
exchangeable measure. We call it the $\nu$-paintbox measure.
\end{enumerate}

Extending seminal results of Kingman \cite{kingman_repre_of_part},
Bertoin has shown in \cite{bertoin_homogeneous} and in his course
at IHP that any exchangeable measure that verifies the required
conditions is a combination of these three types. Hence the
following proposition merely restates these results.

\begin{proposition} \label{prop structures des mesures}
For each exchangeable measure $C$ on $\P$ such that $C(\{\0\})=0,$
and $C(\P \backslash Q(\0,n)) <\infty, \forall n \in \N$ there
exists a unique $\cc \ge 0 $ and a unique measure $\nuc$ on
$\sfleche$ such that
\begin{eqnarray}
\nonumber  \nuc (\{\0 \})= 0, \\
 \label{condition sur nu c} \int_{\sfleche}\left( \sum_{i=1}^{\infty}
x_i^2\right) \nuc (dx) <\infty , \\
\nonumber  \text{ and } C =  \cc \kappa + \mu_{\nuc}.
\end{eqnarray}

For each exchangeable measure $F$ on $\P$ such that $F(\{\1\}) =0$
and $F(\P \backslash Q(\1,n)) <\infty, \forall n \in \N$ there
exists a unique $\ce \ge 0$ and a unique measure $\nuf$ on $\P$
such that
\begin{eqnarray}
\nonumber  \nuf (\{(1,0,..)\})=0, \\
\label{condition sur nu f} \int_{\sfleche} \left(
1-\sum_{i=1}^{\infty} x_i^2\right) \nuf (dx) < \infty , \\
\nonumber \text{and } F = \ce \meros + \mu_{\nuf}.
\end{eqnarray}
\end{proposition}

The two integrability conditions on $\nuf$ and $\nuc$
(\ref{condition sur nu f}) and (\ref{condition sur nu c}) ensure
that $C(\P \backslash Q(\0,n)) <\infty$ and $F(\P \backslash
Q(\1,n)) <\infty.$ See \cite{bertoin_homogeneous} for the
demonstration concerning $F$. The part that concerns $C$ can be
shown by the same arguments.

The condition on $\nuf$ (\ref{condition sur nu f}) may seem at
first sight different from the condition that Bertoin imposes in
\cite{bertoin_homogeneous} and which reads
$$\int_{\sfleche}(1-x_1)\nuf(dx) <\infty$$ but they are in
fact equivalent because $$ 1-\sum_i x_i^2 \le 1 -x_1^2 \le
2(1-x_1)$$ and on the other hand $$1-\sum_i x_i^2  \ge 1-x_1
\sum_i x_i \ge 1-x_1.$$

Thus the above proposition implies that for each EFC process $\Pi$
there is a unique exchangeable fragmentation $\Pi^{(F)}(t)$ and a
unique exchangeable coalescence $\Pi^{(C)}(t)$ such that $\Pi$ is
a combination of $\Pi^{(F)}$ and $\Pi^{(C)}$. This was not obvious
a priori because some kind of compensation phenomena could have
allowed weaker integrability conditions.

One can sum up the preceding analysis in the following
characterization of exchangeable fragmentation-coalescence
processes.

\begin{proposition} \label{existence}
The distribution of an EFC process $\Pi(\cdot)$ is completely
characterized by the initial condition (i.e., the law of
$\Pi(0)$), the measures $\nuf$ and $\nuc$ as above and the
parameters $\ce,\cc \in \R+$.
\end{proposition}

\noindent \textbf{Remark :} The above results are well known for
pure fragmentation or pure coalescence. If, for instance, we
impose $F(\P)=0$ (i.e., there is only coalescence and no
fragmentation, the EFC process is degenerated), the above
proposition shows that our definition agrees with Definition 3 in
Schweinsberg \cite{schweinsberg_multiple}.  On the other hand if
there is only fragmentation and no coalescence, our definition is
equivalent to that given by Bertoin in \cite{bertoin_homogeneous},
which relies on some fundamental properties of the semi-group.
There, the Markov chain property of the restrictions is deduced
from the definition as well as the characterization of the
distribution by $c$ and $\nuf$.

Nevertheless, the formulation of Definition \ref{fundamental
definition} is new. More precisely, it was not known that the
exchangeability and Markov requirement for the restrictions to
$[n]$ was enough to obtain a fragmentation procedure in which each
fragment splits independently from the others (point 2 of
Proposition \ref{prop taux de sauts}).

\subsection{Poissonian construction}

As for exchangeable fragmentation or coalescence, one can
construct EFC processes by using Poisson point processes (PPP in
the following). More precisely let $P_C=((t,\pi^{(C)}(t)),t\ge 0)$
and $P_F=((t,\pi^{(F)}(t),k(t)),t\ge 0)$ be two independent PPP in
the same filtration. The atoms of the PPP $P_C$ are points in $
\R^+ \times \P$ and its intensity measure is given by $dt \otimes
(\mu_{\nuc} + \cc \kappa).$ The atoms of $P_F$ are points in $\R^+
\times \P \times \N $ and its intensity measure  is $dt \otimes
(\ce \meros + \mu_{\nuf}) \otimes \# $ where $\#$ is the counting
measure on $\N$ and $dt$ is the Lebesgue measure.

Take $\pi \in \P$ an exchangeable random  variable and define a
family of $\P_n$-valued processes $\Pi^n(\cdot)$ as follows: for
each $n$ fix $\Pi^n(0)=\pi_{|[n]}$ and
\begin{itemize}
\item if $t$ is not an atom time neither for $P_C$ or $P_F$ then
$\Pi^n(t)=\Pi^n(t-),$ \item if $t$ is an atom time for $P_C$ such
that $(\pi^{(C)}(t))_{|[n]} \neq \0_n$ then $$\Pi^n(t)=
\coag(\Pi^n(t-),\pi^{(C)}(t)),$$ \item if $t$ is an atom time for
$P_{F}$ such that $k(t) <n$ and $(\pi^{(F)}(t))_{|[n]} \neq \1_n$
then
$$\Pi^n(t)= \frag(\Pi^n(t-),\pi^{(F)}(t),k(t)).$$
\end{itemize}

Note that the $\Pi^n$ are well defined because on any finite time
interval, for each $n,$ one only needs to consider a finite number
of atoms. Furthermore $P_C$ and $P_F$ being independent in the
same filtration, almost surely there is no $t$ which is an atom
time for both PPP's. This family is constructed to be compatible
and thus defines uniquely a process $\Pi$ such that $\Pi_{|[n]}
=\Pi^n$ for each $n$. By analogy with exchangeable fragmentations
(\cite{bertoin_homogeneous}) and exchangeable coalescence
(\cite{pitman_multiple,schweinsberg_multiple}) the following
should be clear.

\begin{proposition}\label{construction PPP}
The process $\Pi$ constructed above is an EFC process with
characteristics $\cc,\nuc,\ce$ and $\nuf$.
\end{proposition}

\begin{proof}
It is straightforward to check that the restrictions
$\Pi_{|[n]}(t)$ are Markov chains whose only jumps are either
coagulations or fragmentations. The transition rates are
constructed to correspond to the characteristics $\cc,\nuc,\ce$
and $\nuf$. The only thing left to check is thus exchangeability.
Fix $n \in \N$ and $\sigma$ a permutation of $[n],$ then
$(\sigma(\Pi^n(t)))_{t \ge 0}$ is a jump-hold Markov process. Its
transition rates are given by $q_n^{(\sigma)}(\pi,\pi')
=q_n(\sigma^{-1}(\pi),\sigma^{-1}(\pi')).$

Suppose first that $\pi'=\frag(\pi,\pi'',k)$ for some $\pi''.$
Remark that there exists a unique $l \le \# \pi$ and a permutation
$\sigma'$ of $[m]$ (where $m = |\pi_k|$ is the cardinal of the
$k$-th block of $\pi$ we want to split) such that
$$\sigma^{-1}(\pi')= \frag(\sigma^{-1}(\pi),\sigma'(\pi''),l).$$
Using Proposition \ref{prop taux de sauts} we then obtain that
\begin{eqnarray*}
q_n^{(\sigma)}(\pi,\pi') &=&
q_n(\sigma^{-1}(\pi),\sigma^{-1}(\pi')) \\ &=&
q_n(\sigma^{-1}(\pi),\frag(\sigma^{-1}(\pi),\sigma'(\pi''),l)) \\
&=& F_{m} (\sigma'(\pi'')) \\ &=& F_{m} (\pi'') \\&=&
q_n(\pi,\pi')
\end{eqnarray*}

The same type of arguments show that when $\pi'=\coag(\pi,\pi'')$
for some $\pi''$ we also have
$$q^{(\sigma)}_n(\pi,\pi')=q_n(\pi,\pi').$$ Thus, $\Pi^n$ and
$\sigma(\Pi^n)$ have the same transition rates and hence the same
law.

As this is true for all $n$, it entails that $\Pi$ and
$\sigma(\Pi)$ also have the same law.

\end{proof}

Let $\Pi(\cdot)$ be an EFC process and define $P_t$ its
semi-group, i.e., for a continuous function $\phi : \P \mapsto \R$
$$P_t \phi (\pi) := \mathbf{E}_{\pi} (\phi(\Pi(t)))$$ the
expectation of $\phi(\Pi(t))$ conditionally on $\Pi(0)=\pi.$
\begin{corollary}
An EFC process $\Pi(\cdot)$ has the Feller property, i.e.,
\begin{itemize}
\item for each continuous function $\phi : \P \mapsto \R,$ for
each $\pi \in \P$ one has $$\lim_{t \rightarrow 0+} P_t \phi
(\pi)= \phi(\pi),$$ \item for all $t >0$ the function $\pi \mapsto
P_t \phi (\pi)$ is continuous.
\end{itemize}
\end{corollary}

\begin{proof}
Call $C_f$ the set of functions
$$C_f= \{ f : \P \mapsto \R : \exists n \in \N \text{ s.t. }
\pi_{|[n]} = \pi'_{|[n]} \Rightarrow f(\pi) =f(\pi') \}$$ which is
dense in the space of continuous functions of $\P \mapsto \R$. The
first point is clear for a function $\Phi \in C_f$ (because the
first jump-time of $\Phi(\Pi(\cdot))$ is distributed as an
exponential variable with finite mean). We conclude by density.
For the second point, consider $\pi,\pi' \in \P$ such that
$d(\pi,\pi')< 1/n$ (i.e., $\pi_{|[n]} = \pi'_{|[n]}$) then use the
same PPP $P_C$ and $P_F$ to construct two EFC processes,
$\Pi(\cdot)$ and $\Pi'(\cdot)$, with respective starting points
$\Pi(0)=\pi$ and $\Pi'(0)=\pi'.$ By construction
$\Pi_{|[n]}(\cdot)=\Pi'_{|[n]}(\cdot)$ in the sense of the
identity of the paths.  Hence $$\forall t \ge 0, d(\Pi(t),\Pi'(t))
< 1/n.$$
\end{proof}

Hence, when considering an EFC process, one can always suppose
that one works in the usual augmentation of the natural filtration
$\F_t$ which is then right continuous.

As a direct consequence, one also has the following
characterization of EFC's in terms of the infinitesimal generator
: Let $(\Pi(t),t\ge 0)$ be an EFC process, then the infinitesimal
generator of $\Pi$, denoted by $\mathcal{A}$, acts on the
functions $f \in C_f$ as follows:
\begin{eqnarray*}
\forall \pi \in\P , \mathcal{A} (f)(\pi) = &&\int_{\P} C(d\pi')
(f(\coag(\pi,\pi'))-f(\pi)) \\ &&+ \sum_{k \in \N} \int_{\P}
F(d\pi') (f(\frag(\pi,\pi',k))-f(\pi)),
\end{eqnarray*}
where $F = \ce \meros + \mu_{\nuf}$ and $C = \cc \kappa +
\mu_{\nuc}.$ Indeed, take $f \in C_f$ and $n$ such that
$\pi_{|[n]} = \pi'_{|[n]} \Rightarrow f(\pi) =f(\pi') $, then as
$\Pi_{|[n]}(\cdot)$ is a Markov chain the above formula is just
the usual generator for Markov chains. Transition rates have thus
the required properties and hence this property characterizes EFC
processes.

\subsection{Asymptotic frequencies}

When $A$ is a subset of $\N$ we will write
$$\bar \lambda_A=\limsup_{n \rightarrow \infty} \frac{\#\{ k
\le n : k \in A \}}{n}$$ and $$\underline \lambda_A=\liminf_{n
\rightarrow \infty} \frac{\#\{ k \le n : k \in A \}}{n}.$$ When
the equality $\bar \lambda_A=\underline \lambda_A$ holds we call
$\|A\|,$ the asymptotic frequency of $A$, the common value which
is also the limit
$$\|A\|=\lim_{n \rightarrow \infty} \frac{\#\{ k \le n : k \in
A \}}{n}.$$

If all the blocks of $\pi =(B_1,B_2,..) \in \P$ have an asymptotic
frequency we define
$$\Lambda(\pi)=(\|B_1\|,\|B_2\|,..)^{\downarrow}$$ the decreasing
rearrangement of the $\|B_i\|$'s.

\begin{theorem}\label{Markov pour X}
Let $\Pi(t)$ be an EFC process. Then $$X(t)=\Lambda(\Pi(t))$$
exists almost surely simultaneously for all $t \ge 0,$ and
$(X(t),t\ge 0)$ is a Feller process.
\end{theorem}
The proof (see section 6), which is rather technical, uses the
regularity properties of EFC processes and the existence of
asymptotic frequencies simultaneously for all rational time $t \in
\Q$. We call the process $X(t)$ the associated ranked-mass EFC
process.

\hspace{0,5cm}

\noindent \textbf{Remark} The state space of a ranked mass EFC
process $X$ is $\sfleche.$ Thus, our construction of EFC processes
$\Pi$ in $\P$ started from $\0$ gives us an entrance law
$Q_{(0,0,...)}$ for $X$. More precisely, there is the identity
$Q_{(0,0,...)}(t+s)=Q(s)P_{\0}(t).$ The ranked frequencies of an
EFC process started from $\0$ defines a process with this entrance
law that comes from dust at time $0+$, i.e., the largest mass
vanishes almost surely as $t \searrow 0.$ The construction of this
entrance law is well known for pure coalescence process, see
Pitman \cite{pitman_multiple} for a general treatment, but also
Kingman \cite{kingman_genealogy_large_popu} and
Bolthausen-Sznitman \cite[Corollary 2.3]{sznitman} for particular
cases.

\section{Equilibrium measures}

Consider an EFC process $\Pi$ which is not trivial, i.e.,
$\nuc,\nuf,\ce$ and $\cc$ are not zero simultaneously.

\begin{theorem}
There exists a unique (exchangeable) stationary probability
measure $\rho$ on $\P$ and one has
$$\rho=\delta_{\0} \Leftrightarrow \cc=0 \text{ and }
\nuc \equiv 0$$ and
$$\rho=\delta_{\1} \Leftrightarrow \ce=0 \text{ and }
\nuf \equiv 0$$ where $\delta_{\pi}$ is the Dirac mass at $\pi$.

Furthermore, $\Pi(\cdot)$ converges in distribution to $\rho$.

\end{theorem}

\begin{proof}

If the process $\Pi$ is a pure coalescence process (i.e., $\ce=0
\text{ and } \nuf(\cdot) \equiv 0$) it is clear that $\1$ is an
absorbing state towards which the process converges almost surely.
In the pure fragmentation case it is $\0$ that is absorbing and
attracting.

In the non-degenerated case, for each $n \in \N,$ the process
$\Pi_{|[n]}(\cdot)$ is a finite state Markov chain. Let us now
check the irreducibility in the non-degenerated case. Suppose
first that $\nuf(\sfleche)>0.$ For every state $\pi \in \P_n,$ if
$\Pi_{|[n]}(t)=\pi$ there is a positive probability that the next
jump of $\Pi_{|[n]}(t)$ is a coalescence. Hence, for every
starting point $\Pi_{|[n]}(0)=\pi \in \P_n$ there is a positive
probability that $\Pi_{|[n]}(\cdot)$ reaches $\1_n$ in finite time
$T$ before any fragmentation has occurred. Now take $x \in
\sfleche$ such that $x_2>0$ and recall that $\mu_x$  is the
$x$-paintbox distribution. Then for every $\pi \in \P_n$ with $\#
\pi=2$ (recall that $\# \pi$ is the number of non-empty blocks of
$\pi$) one has
$$\mu_x(Q(\pi,n))>0.$$ That is the $n$-restriction of the
$x$-paintbox partition can be any partition of $[n]$ in two blocks
with positive probability. More precisely if $\pi \in \P_n$ is
such that $\pi=(B_1,B_2,\o,\o...)$ with $|B_1|=k$ and $|B_2|=n-k$
then
$$\mu_x(Q(\pi,n)) \ge x_1^k x_2^{n-k} + x_2^k x_1^{n-k}.$$ Hence,
for any $\pi \in \P$ with $\# \pi =2$, the first transition after
$T$ is $\1_n \rightarrow \pi$ with positive probability. As any
$\pi \in \P_n$ can be obtained from $\1_n$ by a finite series of
binary fragmentations we can iterate the above idea to see that
with positive probability the jumps that follow $T$ are exactly
the sequence of binary splitting needed to get to $\pi$ and the
chain is hence irreducible.

Suppose now that $\nuf \equiv 0$, there is only erosion $\ce >0,$
and that at least one of the following two condition holds
\begin{itemize}
\item for every $k \in \N$ one has $\nuc(\{ x \in \sfleche :
\sum_{i=1}^{i=k} x_i <1 \})>0$, \item there is a Kingman
component, $\cc >0$,
\end{itemize}
then almost the same demonstration applies. We first show that the
state $\0_n$ can be reached from any starting point by a series of
splittings corresponding to erosion, and that from there any $\pi
\in \P_n$ is reachable through binary coagulations.

In the remaining case (i.e., $\cc=0,\nuf\equiv 0$ and there exists
$k >0$ such that $\nuc(\{ x \in \sfleche : \sum_{i=1}^{i=k} x_i <1
\})=0$) the situation is slightly different in that $\P_n$ is not
the irreducible class. It is easily seen that the only partitions
reachable from $\0_n$ are those with at most $k$ non-singletons
blocks. But for every starting point $\pi$ one reaches this class
in finite time almost surely. Hence there is no issues with the
existence of an invariant measure for this type of $\Pi_{|[n]},$
it just does not charge partitions outside this class.

Thus there exists a unique stationary probability measure
$\rho^{(n)}$ on $\P_n$ for the process $\Pi_{|[n]}.$ Clearly by
compatibility of the $\Pi_{|[n]}(\cdot)$ one must have
$$\text{Proj}_{\P_n}(\rho^{(n+1)})(\cdot) = \rho^{(n)}(\cdot)$$ where
$\text{Proj}_{\P_n}(\rho^{(n+1)})$ is the image of $\rho^{(n+1)}$
by the projection on $\P_n.$ This implies that there exists a
unique probability measure $\rho$ on $\P$ such that for each $n$
one has $\rho^{(n)}(\cdot)=\text{Proj}_{\P_n}(\rho)(\cdot).$ The
exchangeability of $\rho$ is a simple consequence of the
exchangeability of $\Pi$. Finally, the chain $\Pi_{|[2]}(\cdot)$
is specified by two transition rates $\{1\}\{2\} \rightarrow
\{1,2\}$ and $\{1,2\} \rightarrow \{1\}\{2\},$ which are both
non-zero as soon as the EFC is non-degenerated. Hence,
$$\text{Proj}_{\P_2}(\rho) (\cdot) \not \in \{ \delta_{\1_2}(\cdot) , \delta_{\0_2}(\cdot)\}.$$
Hence, when we have both coalescence and fragmentation $\rho \not
\in \{\delta_{\1},\delta_{\0}\}.$

The $\Pi_{|[n]}(\cdot)$ being finite states Markov chains, it is
well known that they converge in distribution to $\rho^{(n)},$
independently of the initial state. By definition of the
distribution of $\Pi$ this implies that $\Pi(\cdot)$ converges in
distribution to $\rho.$

\end{proof}

Although we cannot give an explicit expression for $\rho$ in terms
of $\cc, \nuc,$ $\ce$ and $\nuf,$ we now relate certain properties
of $\rho$ to these parameters. In particular we will ask ourselves
the following two natural questions:
\begin{itemize}
\item under what conditions does $\rho$ charge only partitions
with an infinite number of blocks, resp. a finite number of
blocks, resp. both ? \item under what conditions does $\rho$
charge partitions with dust (i.e., partitions such that $\sum_i
\|B_i\| \le 1$ where $\|B_i\|$ is the asymptotic frequency of
block $B_i$) ?
\end{itemize}

The proofs of the results in the remaining of this section are
placed in section 6.

\subsection{Number of blocks}

We will say that an EFC process \textit{fragmentates quickly} if
$\ce>0$ or $\nuf(\sfleche) =\infty.$ If it is not the case (i.e.,
$\ce=0$ and $\nuf(\sfleche) <\infty$) we say that it
\textit{fragmentates slowly.}

We first examine whether of not $\rho$ charges partitions with a
finite number of blocks.

\begin{theorem}\label{propo number}
\begin{enumerate}

\item Let $\Pi(\cdot)$ be an EFC process that fragmentates
quickly. Then
$$\rho (\{ \pi \in \P : \# \pi <\infty\}) = 0.$$

\item Let $\Pi(\cdot)$ be an EFC process that fragmentates slowly
and such that
$$\nuf(\{x \in \sfleche : x_1+x_2<1\})=0$$ (the fragmentation
component is binary),  then
$$\cc >0 \Rightarrow \rho(\{\pi \in \P : \#\pi< \infty\}) =1.$$
\end{enumerate}
\end{theorem}

Thus for an EFC process with a binary fragmentation component, a
Kingman coalescence component and no erosion (i.e., $ \cc>0,\ce
=0$ and $\nuf(\{x \in \sfleche : x_1+x_2<1\})=0$) we have the
equivalence
$$\rho(\{\pi \in \P : \# \pi =\infty\}) =1 \Leftrightarrow
\nuf(\sfleche)=\infty$$ and when $\nuf(\sfleche)<\infty$ then
$\rho(\{\pi \in \P : \# \pi =\infty\}) =0.$

\subsection{Dust}

For any fixed time $t$ the partition $\Pi(t)$ is exchangeable.
Hence, by Kingman's theory of exchangeable partition, its law is a
mixture of paintbox processes. A direct consequence is that every
block $B_i(t)$ of $\Pi(t)$ is either a singleton or an infinite
block with strictly positive asymptotic frequency. Recall that the
asymptotic frequency of a block $B_i(t)$ is given by
$$\|B_i(t)\|=\lim_{n \rightarrow \infty} \frac{1}{n} \text{Card}\{ k
\le n : k \in B_i(t) \}$$ so part of Kingman's result is that this
limit exists almost surely for all $i$ simultaneously. The
asymptotic frequency of a block corresponds to its mass, thus
singletons have zero mass, they form what we call dust. More
precisely, for $\pi \in \P$ define the set
\begin{eqnarray*}
\text{dust}(\pi)&:=& \bigcup_{j : \|B_j\|=0}B_j .
\end{eqnarray*}
When $\pi$ is exchangeable we have almost surely $$ \text{dust}
(\pi)=  \{i \in \N : \exists j \text{ s.t. } \{i\} = B_j\} $$ and
$$\sum_i \|B_i\| + \|\text{dust}(\pi)\|=1.$$

For fragmentation or EFC processes, dust can be created
\textit{via} two mechanisms: either from erosion (that's the atoms
that correspond to the erosion measure $\ce \meros$ when $\ce
>0$), or from sudden splitting which corresponds to atoms
associated to the measure $\mu_{\nuf'}$ where $\nuf'$ is simply
$\nuf$ restricted to $\{ s \in \sfleche : \sum_i s_i <1\} $.
Conversely, in the coalescence context mass can condensate out of
dust, thus giving an entrance law in $\sfleche$, see
\cite{pitman_multiple}.

The following theorem states that when the coalescence is strong
enough in an EFC process, the equilibrium measure does not charge
partitions with dust. We say that an EFC process \textit{coalesces
quickly} (resp. \textit{slowly}) if $ \int_{\sfleche} \left(
\sum_i x_i \right) \nuc(dx) =\infty \text{ or } \cc>0$ (resp. $
\int_{\sfleche} \left( \sum_i x_i \right) \nuc(dx) <\infty$ and
$\cc=0$).

\begin{theorem}\label{theorem dust equ}
Let $(\Pi(t),t\ge0)$ be an EFC process that coalesces quickly and
$\rho$ its invariant probability measure. Then $$\rho(\{\pi\in \P
: \text{dust}(\pi) \neq \o\})=0. $$
\end{theorem}
In case of no fragmentation, this follows from Proposition 30 in
\cite{schweinsberg_multiple}.

\subsection{Equilibrium measure for the ranked mass EFC process}
For $\rho$ the equilibrium measure of some EFC process with
characteristics $\nuf,\nuc,\ce$ and $\cc,$ let $\theta$ be the
image of $\rho$ by the map $\P \mapsto \sfleche : \pi \mapsto
\Lambda(\pi)$.

\begin{proposition}
Let $X$ be a ranked-mass EFC process with characteristics
$\nuf,\nuc,\ce$ and $\cc.$ Then $\theta$ is its unique invariant
probability measure.
\end{proposition}
\begin{proof}
As for each fixed $t$ one has $$P_{\rho} (\Lambda(\Pi(t)) \in A) =
\rho(\{ \pi : \Lambda(\pi) \in A  \} )=\theta(A)$$ it is clear
that $\theta$ is an invariant probability measure.

Suppose that $\theta$ is an invariant measure for $X$ and fix $t
\ge 0.$ Hence if $X(0)$ has distribution $\theta$ so does
$X(t)=\Lambda(\Pi(t)).$ As $\Pi(t)$ is exchangeable it is known by
Kingman's theory of exchangeable partitions that $\Pi(t)$ has law
$\mu_{\theta}(\cdot)$ the mixture of paintbox processes directed
by $\theta.$ This implies that $\mu_{\theta}(\cdot)$ is invariant
for $\Pi$ and hence $\mu_{\theta}(\cdot)=\rho(\cdot)$ and thus
$\theta$ is the unique invariant measure for $X.$
\end{proof}

\section{Path properties}

\subsection{Number of blocks along the path.}
One of the problem tackled by Pitman \cite{pitman_multiple} and
Schweinsberg \cite{schweinsberg_comes_down,schweinsberg_multiple}
about coalescent processes is whether or not they come down from
infinity. Let us first recall some of their results. By definition
if $\Pi^C(\cdot)$ is a standard coalescent $\Pi^C(0)=\0$ and thus
$\# \Pi^C(0)=\infty.$ We say that $\Pi^C$ comes down from infinity
if $\# \Pi^C(t) <\infty$ a.s. for all $t>0$. We say it stays
infinite if $\# \Pi^C(t) =\infty$ a.s. for all $t>0$.

Define $\Delta_f :=\{x \in \sfleche : \exists i \in \N \text{ s.t.
} \sum_{j =1}^i x_j =1\}$. We know by Lemma 31 in
\cite{schweinsberg_comes_down}, which is a generalization of
Proposition 23 in \cite{pitman_multiple}, that if
$\nuc(\Delta_f)=0$ the coalescent either stays infinite or comes
down from infinity.

For $b \ge 2$ let $\lambda_b$ denote the total rate of all
collisions when the coalescent has $b$ blocks $$\lambda_b=
\mu_{\nuc} (\P \backslash Q(\0,b))+\cc \frac{b(b-1)}{2}.$$ Let
$\gamma_b$ be the total rate at which the number of blocks is
decreasing when the coalescent has $b$ blocks, $$ \gamma_b =
\cc\frac{b(b-1)}{2} + \sum_{k=1}^{b-1} (b-k) \mu_{\nuc}(\{\pi : \#
\pi_{|[b]} = k\}).$$

If $\nuc (\Delta_f) = \infty$ or $\sum_{b =2}^{\infty}
\gamma_b^{-1} < \infty,$ then the coalescent comes down from
infinity. The converse is not always true but holds for instance
for the important case of the $\Lambda$-coalescents (i.e., those
for which many fragments can merge into a single block, but only
one such merger can occur simultaneously).

This type of properties concerns the paths of the processes, and
it seems that they bear no simple relations with properties of the
equilibrium measure. For instance the equilibrium measure of a
coalescent that stays infinite is $\delta_{\1}(\cdot)$ and
therefore only charges partitions with one block, but its path
lays entirely in the subspace of $\P$ of partitions with an
infinite number of blocks.

\vspace{1cm}

Let $\Pi(\cdot)$ be an EFC process. Define the sets $$G:=\{t \ge 0
: \#\Pi(t)=\infty \}$$ and  $$\forall k \in \N, G_k :=\{t \ge 0 :
\#\Pi(t) > k  \}.$$ Clearly every arrival time $t$ of an atom of
$P_C$ such that $\pi^{(C)}(t) \in \Delta_f$ is in $G^c$ the
complementary of $G$. In the same way an arrival time $t$ of an
atom of $P_F$ such that $\pi^{(F)}(t) \in \sfleche \backslash
\Delta_f$ and $B_{k(t)}(t-)$ (the fragmented block) is infinite
immediately before the fragmentation, must be in $G$. Hence, if
$\nuf(\sfleche \backslash \Delta_f) =\infty$ and $\nuc(\Delta_f)
=\infty,$  then both $G$ and $G^c$ are everywhere dense, and this
independently of the starting point which may be $\1$ or $\0$.

The following proposition shows that when the fragmentation rate
is infinite, $G$ is everywhere dense. Recall the notation
$\Pi(t)=(B_1(t),B_2,(t),...)$.

\begin{theorem}
Let $\Pi$ be an EFC process that fragmentates quickly. Then, a.s.
$G$ is everywhere dense.
\end{theorem}

As $G = \cap G_k$ we only need to show that a.s. for each $k \in
\N$ the set $G_k$ is everywhere dense and open to conclude with
Baire theorem. The proof relies on two lemmas.

\begin{lemma}
Let $\Pi$ be an EFC process that fragmentates quickly started from
$\1.$  Then, a.s. for all $k \in \N$
$$\inf \{ t\ge 0 : \# \Pi(t) >k\} =0.$$
\end{lemma}
\begin{proof}

Fix $k \in \N$ and $\epsilon >0$, we are going to show that there
exists $t \in \left[0,\epsilon\right[$ such that $$\exists n\in \N
: \# \Pi_{|[n]}(t) \ge k.$$

Recall the notation $B(i,t)$ for the block of $\Pi(t)$ that
contains $i.$ As $\nuf(\sfleche) =\infty$ (or $\ce >0$) it is
clear that almost surely $\exists n_1 \in \N : \exists t_1 \in
\left[0,\epsilon\right[$ such that
$\Pi_{|[n_1]}(t_1-)=\1_{|[n_1]}$ and $t_1$ is a fragmentation time
such that $\Pi_{|[n_1]}(t_1)$ contains at least two blocks, say
$B(i_1,t_1)\cap [n_1]$ and $B(i_2,t_1)\cap [n_1],$ of which at
least one is not a singleton and is thus in fact infinite when
seen in $\N$. The time of coalescence of $i_1$ and $i_2$ (i.e.,
the first time at which they are in the same block again is
exponentially distributed with parameter $$ \int_{\sfleche}
(\sum_i x_i^2) \nuc (dx) + \cc <\infty.$$ Hence if we define
$$\tau_{i_1,i_2}(t_1):=\inf \{ t \ge t_1 : i_1
\overset{\Pi(t)}{\sim} i_2 \}$$ then almost surely we can find
$n_2>n_1$ large enough such that the first time $t_2$ of
fragmentation of $B(i_1,t_1)\cap [n_2]$ or $B(i_2,t_1)\cap [n_2]$
is smaller than $\tau_{i_1,i_2}(t_1)$ (i.e., $i_1$ and $i_2$ have
not coalesced yet) and $t_2$ is a fragmentation time at which
$B(i_1,t_2-)\cap [n_2]$ or $B(i_2,t_2-)\cap [n_2]$ is split into
two blocks. Hence at $t_2$ there are at least $3$ non-empty blocks
in $\Pi_{|[n_2]}(t_2)$, and at least one of them is not a
singleton. By iteration, almost surely, $\exists n_k : \exists t_k
\in \left[0,\epsilon\right[$ such that $t_k$ is a fragmentation
time and $$\# \Pi_{|[n_k]}(t_k) \ge k.$$

\end{proof}

\begin{lemma}\label{lemma gk dense}
Let $\Pi$ be an EFC process that fragmentates quickly. Then, a.s.
$G_k$ is everywhere dense and open for each $k \in \N$.
\end{lemma}
\begin{proof}
Fix $k \in \N$, call $\Gamma_k=\{t^{(k)}_1 < t^{(k)}_2 <...\}$ the
collection of atom times of $P_C$ such that a coalescence occurs
on the $k+1$ first blocks if there are more than $k+1$ blocks,
i.e.,$$\pi^{(C)}(t) \not \in Q(\0,k+1)$$ (recall that
$Q(\0,k+1)=\{ \pi \in \P : \pi_{|[k+1]}=\0_{k+1}\}$). Suppose $t
\in G_k$, then by construction $\inf \{ s > t : s \in G_k^c\} \in
\Gamma_k$ (because one must at least coalesce the first $k+1$
distinct blocks present at time $t$ before having less than $k$
blocks.  As the $t^{(k)}_i$ are stopping times, the strong Markov
property and the first lemma imply that $G_k^c \subseteq
\Gamma_k.$ Hence $G_k$ is a dense open subset of $\R+$.
\end{proof}

We can apply Baire's theorem to conclude that $\cap_k G_k=G$ is
almost surely everywhere dense in $\R+$.

As a corollary, we see that when the coalescence is ``mostly"
Kingman (i.e., $\cc >0, \nuc(\sfleche) <\infty$) and the process
fragmentates quickly ($\ce>0$ or $\nuf(\sfleche)=\infty$), then we
have that the set of times $G^c$ is exactly the set of atom times
for $P_C$ such that $\pi^{(C)}(\cdot) \in \Delta_f$. Define
$\Delta_f(k):=\{x \in \Delta_f : \sum_1^k x_i =1\}.$
\begin{corollary}\label{corollaire}
Consider an EFC process that fragmentates quickly. When $\cc \ge
0$ and $\nuc(\sfleche)<\infty$ one has
\begin{eqnarray} \label{ensemble G}
G^c= \{ t : \pi^{(C)}(t) \in \Delta_f\}
\end{eqnarray}
and for all $n \ge 2$
\begin{eqnarray} \label{ensemble Gn}
G^c_n=\{ t : \pi^{(C)}(t) \in \Delta_f(n-1)\}.
\end{eqnarray}
\end{corollary}

\begin{proof}
As $G^c=\cup G^c_n,$ it suffices to show (\ref{ensemble Gn}) for
some $n \in \N.$

Recall from the proof of Lemma \ref{lemma gk dense} that $G^c_k
\subseteq \Gamma_k=\{t_1^{(k)},t_2^{(k)},..\}$ the set of
coalescence times at which $\pi^{(C)}(t) \not \in Q(\0,k+1).$

Now fix $n \in \N$ and consider simultaneously the sequence
$(t_i^{(k)})_{i \in \N}$ and $(t_i^{(k+n)})_{i \in \N}.$ It is
clear that for each $i \in \N, \exists j \in \N$ such that $$
t_i^{(k)} = t_j^{(k+n)}$$ because $\pi^{(c)}(t) \not \in Q(\0,k+1)
\Rightarrow \pi^{(c)}(t) \not \in Q(\0,k+n+1).$ Furthermore the
$t_i^{(k+n)}$ have no other accumulation points than $\infty,$
thus there exists $r_1^{(k+n)}  < t_1^{(k)}$ and $n_1 < \infty $
such that for all $s \in ]r_1^{(k+n)},t_1^{(k)}[ :
\#\Pi_{[n_1]}(s) > k+n.$ Hence, a necessary condition to have $\#
\Pi_{|[n_1]}(t_1^{(k)}) < k $ is that $t_1^{(k)}$ is a multiple
collision time, and more precisely $t_1^{(k)}$ must be a collision
time such that $\# \pi^{(c)}_{|[k+n]}(t_1^{(k)}) \le k.$ Hence
$$G_k^c \subseteq \{ t_i^{(k)} \text{ s.t. } \# \pi^{(c)}_{|[k+n]}(t_i^{(k)}) \le
k\}.$$As this is true for each $n$ almost surely, the conclusion
follows.
\end{proof}

\hspace{1cm}

As recently noted by Lambert \cite{amaury}, there is an
interpretation of some EFC processes in terms of population
dynamics. More precisely if we consider an EFC process $(\Pi(t),t
\ge 0)$ such that $\nuf(\sfleche)<\infty$ and
\begin{equation*}
(H) \qquad \qquad \left\{ \begin{array}{c}  \nuf(\sfleche
\backslash \Delta_f)=0
\\\ce=0 \\\nuc(\sfleche)=0 \\ \cc >0
\end{array} \right.
\end{equation*}
then, if at all time all the blocks of $\Pi(t)$ are infinite we
can see the number of blocks $(Z(t)=\#\Pi(t),t\ge 0)$ as the size
of a population where each individuals gives rise (without dying)
to a progeny of size $i$ with rate $\nuf(\Delta_f(i+1))$ and there
is a negative density-dependence due to competition pressure. This
is reflected by the Kingman coalescence phenomena which results in
a quadratic death rate term. The natural death rate is set to 0,
i.e., there is no linear component in the death rate. In this
context, an EFC process that comes down from infinity corresponds
to a population started with a very large size. Lambert has shown
that a sufficient condition to be able to define what he terms a
logistic branching process started from infinity is $$(L) \qquad
\qquad \sum_k p_k \log k <\infty$$ where $p_k =
\nuf(\Delta_f(k+1)).$

More precisely, this means that if $P_n$ is the law of the
$\N$-valued Markov chain $(Y(t),t\ge 0)$ started from $Y(0)=n$
with transition rates
$$\forall i \in \N \left\{ \begin{array}{c} i \rightarrow i+j \text{ with rate }
i p_j \text{ for all }j \in \N  \\ i \rightarrow i-1 \text{ with
rate } \cc i(i-1)/2 \text{ when }i>1. \end{array}\right. ,$$ then
$P_n$ converge weakly to a law $P_{\infty}$ which is the law of a
$\N \cup \infty$-valued Markov process $(Z(t),t\ge 0) $ started
from $\infty,$ with same transition semi-group on $\N$ as $Y$ and
whose entrance law can be exhibited. Moreover, if we call
$\tau=\inf \{t \ge  0 : Z(t)=1  \}$ we have that
$\E(\tau)<\infty.$

As $\#\Pi(\cdot)$ has the same transition rates as $Y(\cdot)$ and
the entrance law from $\infty$ is unique, these processes have the
same law. Hence the following is a simple corollary of Lambert's
result.
\begin{proposition}
Let $\Pi$ be an EFC process started from dust (i.e., $\Pi(0)=\0$)
and verifying the conditions (H) and (L). Then one has
\begin{eqnarray*}
\label{nuf fini} \forall t>0\text{, }  \# \Pi(t) <\infty  \text{
a.s.}
\end{eqnarray*}
\end{proposition}

\begin{proof}
If $T=\inf\{ t : \#\Pi(t)<\infty\}$, Lambert's result implies that
$\E(T)<\infty$ and hence $T$ is almost surely finite. A simple
application of Proposition 23 in \cite{pitman_multiple} and Lemma
31 in \cite{schweinsberg_multiple} shows that if there exists $t
<\infty$ such that $\#\Pi(t)<\infty$ then $\inf\{ t :
\#\Pi(t)<\infty\}=0.$ To conclude, it is not hard to see that if
$\Pi(0)=\1$ then $\inf \{t \ge 0 : \#\Pi(t)=\infty\}=\infty.$ This
entails that when an EFC process verifying (H) and (L) reaches a
finite level it cannot go back to infinity. As $\inf\{ t :
\#\Pi(t)<\infty\}=0,$ this means that
$$\forall t >0, \#\Pi(t)<\infty.$$
\end{proof}

\textbf{Remark :} Let $\Pi(\cdot)=(B_1(\cdot),B_2(\cdot),...)$ be
a ``(H)-(L)" EFC process started from dust, $\Pi(0)=\0.$ Then for
all $t > 0$ one has a.s. $\sum_i \|B_i(t)\| =1.$ This is clear
because at all time $t>0$ there are only a finite number of
blocks.

If we drop the hypothesis $\nuf(\sfleche)<\infty$ (i.e., we drop
(L) and we suppose $\nuf(\sfleche)=\infty$), the process $\Pi$
stays infinite (Corollary \ref{corollaire}). We now show that
nevertheless, for a fixed $t,$ almost surely $\|B_1(t)\|>0$. We
define by induction a sequence of integers $(n_i)_{i \in \N}$ as
follows: we fix $n_1=1, t_1=0$ and for each $i>1$ we chose $n_i$
such that there exists a time $t_i <t$ such that $t_i$ is a
coalescence time at which the block $1$ coalesces with the block
$n_i$ and such that $n_i >w_{n_{i-1}}(t_{i-1})$ where $w_k(t)$ is
the least element of the $k$th block at time $t$. This last
condition ensures that $(w_{n_i}(t_i))$ is a strictly increasing
sequence because one always has $w_n(t)\ge n.$ The existence of
such a construction is assured by the condition $\cc >0.$ Hence at
time $t$ one knows that for each $i$ there has been a coalescence
between $1$ and $w_{n_i}(t_i)$. Consider $(\Pi^{(F)}_t(s),s\in
[0,t[)$ a coupled fragmentation process defined as follows:
$\Pi^{(F)}_t(0)$ has only one block which is not a singleton which
is
$$B_1^{(F)}(0)=\{1,w_{n_2}(t_2),w_{n_3}(t_3),.....\}.$$ The
fragmentations are given by the same PPP $P_F$ used to construct
$\Pi$ (and hence the processes are coupled). It should be clear
that if $w_{n_i}(t_i)$ is in the same block with $1$ for
$\Pi^{(F)}(t)$ the same is true for $\Pi(t)$ because it means that
no dislocation separates $1$ from $w_{n_i}(t_i)$ during $[0,t]$
for $\Pi^{(F)}$ and hence
$$1 \overset{\Pi(t)}{\sim} w_{n_i}(t_i).$$ Using this fact and standard
properties of homogeneous fragmentations one has a.s.
$$\|B_1(t)\| \ge \|B_1^{(F)}(t)\|>0.$$

Hence for all $t >0$ one has $P(\{1\} \subset
\text{dust}(\Pi(t)))=0$ and hence $P(\text{dust}(\Pi(t)) \neq
\o)=0.$ Otherwise said, when $\nuf(\sfleche) =\infty$ the
fragmentation part does not let a ``(H)" EFC process come down
from infinity, but it let the dust condensates into mass. Note
that ``binary-binary"\footnote{i.e., $\ce = 0, \nuf(\{ x  :
x_1+x_2 <1\}=0, \nuc \equiv 0$ and $\cc>0$.} EFC processes are a
particular case. The question of the case $\nuf(\sfleche)<\infty$
but $(L)$ is not true remains open.

\subsection{Missing mass trajectory}

This last remark prompts us to study more generally the behavior
of the process of the missing mass $$D(t)=
\|\text{dust}(t)\|=1-\sum_i\|B_i(t)\|.$$

In \cite{pitman_multiple} it was shown (Proposition 26) that for a
pure $\Lambda$-coalescence started from $\0$ (i.e., such that
$\nuc(\{ x \in \sfleche : x_2>0\})=0$) $$\xi(t):=-\log (D(t))$$
has the following behavior:
\begin{itemize}
\item either the coalescence is quick ($\cc>0$ or $\int_{\sfleche}
(\sum_i x_i) \nuc (dx) = \infty$) and then $D(t)$ almost surely
jumps from 1 to 0 immediately (i.e., $D(t)=0$ for all $t
>0,$)

\item either the coalescence is slow ($\cc=0$ and $\int_{\sfleche}
(\sum_i x_i) \nuc (dx) <\infty$) and one has that $\xi(t)$ is a
drift-free subordinator whose Lévy measure is the image of
$\nuc(dx)$ via the map $x \mapsto -\log (1-x_1).$
\end{itemize}

In the following we make the following hypothesis about the EFC
process we consider
\begin{equation*}
\text{(H')} \qquad \left\{ \begin{array}{c}  \cc=0
\\ \int_{\sfleche} (\sum_i x_i) \nuc (dx) <\infty \\ \nuf(\{ x
\in \sfleche  : \sum_i x_i <1 \}=0 \
\end{array}\right.
\end{equation*}
The last assumption means that sudden dislocations do not create
dust.

Before going any further we should also remark   that without loss
of generality we can slightly modify the PPP construction given in
Proposition \ref{construction PPP} : We now suppose that $P_F$ is
the sum of two point processes $P_F =P_{Disl} + P_e$ where
$P_{\text{Disl}}$ has measure intensity $\mu_{\nuf} \otimes \#$
and $P_e$ has measure intensity $\ce \otimes \#.$ If $t$ is an
atom time for  $P_{Disl}$ one obtains $\Pi(t)$ from $\Pi(t-)$ as
before, if $P_e$ has an atom at time $t,$ say $(t,k(t)),$ then
$\Pi(t-)$ is left unchanged except for $k(t)$ which becomes a
singleton if this was not already the case. Furthermore, if $t$ is
an atom time for $P_C$ we will coalesce $B_i(t-)$ and $B_j(t-)$ at
time $t$ if and only if $w_i(t-)$ and $w_j(t-)$ (i.e., the least
elements of $B_i(t-)$ and $B_j(t-)$ respectively) are in the same
block of $\pi^{(C)}(t).$ This is equivalent to say that from the
point of view of coalescence the labelling of the block is the
following: if $i$ is not the least element of its block $B_i$ is
empty, and if it is the least element of its block then $B_i$ is
this block. To check this, one can for instance verify that the
transition rates of the restrictions $\Pi_{|[n]}(\cdot)$ are left
unchanged.

\begin{proposition}
Let $\Pi$ be an EFC process verifying (H'). Then $\xi$ is solution
of the SDE
$$d \xi(t):= d\sigma(t) - \ce (e^{\xi(t)}-1)dt$$ where $\sigma$ is a drift-free subordinator whose
Lévy measure is the image of $\nuc(dx)$ via the map $x \mapsto
-\log (1-\sum_i x_i).$

\end{proposition}

The case when $\ce=0$ is essentially a simple extension of
Proposition 26 in \cite{pitman_multiple} which can be shown with
the same arguments. More precisely, we use a coupling argument. If
we call $(\Pi^{(C)}(t),t\ge 0 )$ the coalescence process started
from $\Pi(0)$ and constructed with the PPP $P_C,$ we claim that
for all $t$ $$\text{dust} (\Pi(t)) = \text{dust}(\Pi^{(C)}(t)).$$
This is clear by remarking that for a given $i$ if we define
$$T^{(C)}_i=\inf \{t >0 : i \not \in \text{dust}(\Pi^{(C)}(t))\}$$
we have that $T^{(C)}_i$ is necessarily a collision time which
involves $\{i\}$ and the new labelling convention implies that
$$T^{(C)}_i=\inf \{t >0 : i \not \in \text{dust}(\Pi(t))\}.$$
Furthermore, given a time $t$, if $i \not \in \text{dust}(\Pi(t))$
then $ \forall s \ge 0 : i \not \in \text{dust}(\Pi(t+s)).$ Hence
for all $t \ge 0$ one has
$\text{dust}(\Pi(t))=\text{dust}(\Pi^{(C)}(t))$ and thus
Proposition 26 of \cite{pitman_multiple} applies.

We now concentrate on the case $\ce>0$. Define $$D_n(t) :=
\frac{1}{n} \# \{ \text{dust}(\Pi(t)) \cap [n]\}.$$ Note that
$\text{dust}(\Pi(t)) \cap [n]$ can be strictly included in the set
of the singletons of the partition $\Pi_{|[n]}(t).$ Remark that
the process $D_n$ is a Markov chain with state-space
$\{0,1/n,...,(n-1)/n,1\}.$ We already know that $D$ is a càdlàg
process and that almost surely, for all $t \ge 0$ one has $D_n(t)
\rightarrow D(t).$

First we show that
\begin{lemma}
With the above notations $D_n \Rightarrow D.$
\end{lemma}
\begin{proof}
One only has to show that the sequence $D_n$ is tight because we
have convergence of the finite dimensional marginal laws (see for
instance \cite[VI.3.20]{jacod_shiryaev}).

The idea is to use Aldous' tightness criterion
(\cite[VI.4.5]{jacod_shiryaev}). The processes $D_n$ are bounded
by $0$ and $1$ and hence the first condition is trivial. We have
to check that $\forall \epsilon >0$ $$\lim_{\theta \searrow 0}
\limsup_n \sup_{S,T \in \mathcal T_N^n ; S \le T \le S +\theta}
P(|D_n(T)-D_n(S)| \ge \epsilon) =0$$ where $\mathcal T_N^n$ is the
set of all stopping times in the natural filtration of $D_n$
bounded by $N.$

First remark that
\begin{eqnarray*}
\sup_{S,T \in \mathcal T_N^n ; S \le T \le S +\theta} &&
P(|D_n(T)-D_n(S)| \ge \epsilon) \\ &\le&  \sup_{S \in \mathcal
T_N^n } P( \sup_{t \le \theta }|D_n(S+t)-D_n(S)| \ge \epsilon)
\end{eqnarray*}
hence we will work on the right hand term.

Fix $S \in \mathcal T_N^n.$ First we wish to control $P( \sup_{t
\le \theta }(D_n(S+t)-D_n(S)) \ge \epsilon).$ Remark that the
times $t$ at which $\Delta(D_n(t))=D_n(t)-D_n(t-)>0$ all are atom
times of $P_F$ such that $\pi^{(F)}(t)=\epsilon_i$ for some $i \le
n$ (recall that $\epsilon_i$ is the partition of $\N$ that
consists in two blocks: $\{i\}$ and $\N\backslash \{i\}$). Hence,
clearly,
\begin{eqnarray*}
P( \sup_{t \le \theta }(D_n(S+t)-D_n(S)) \ge \epsilon) &\le&
P(\frac{1}{n}\sum_{s \in [S,S+\theta]}
\indic{\pi^{(F)}(s)=\epsilon_i , i=1,...,n} \ge \epsilon) .
\end{eqnarray*}
The process
$$\Big(\sum_{i=1}^n \sum_{s\in [S,S+\theta]}
\indic{\pi^{(F)}(s)=\epsilon_i}\Big)_{\theta \ge 0}$$ is a sum of
$n$ independent standard Poisson processes with intensity $\ce,$
hence for each $\eta>0$ and $\epsilon >0$ there exists $\theta_0$
and $n_0$ such that for each $\theta \le \theta_0$ and $n\ge n_0$
one has
$$P(\frac{1}{n}  \sum_{i=1}^n \sum_{s\in [S,S+\theta]} \indic{\pi^{(F)}(s)=\epsilon_i}
>\epsilon)= P(\frac{1}{n}  \sum_{i=1}^n \sum_{s\in [0,0+\theta]} \indic{\pi^{(F)}(s)=\epsilon_i}
>\epsilon) <\eta$$ where the first equality is just the strong
Markov property in $S.$ Hence, the bound is uniform in $S$ and one
has that for each $\theta \le \theta_0$ and $n\ge n_0$
$$\sup_{S \in \mathcal
T_N^n } P( \sup_{t \le \theta }(D_n(S+t)-D_n(S)) \ge \epsilon)
<\eta.$$

Let us now take care of $P( \sup_{t \le \theta }(D_n(S)-D_n(S+t))
\ge \epsilon).$ We begin by defining a coupled coalescence process
as follows: we let $\Pi_S^{(C)}(0)= \0,$ and the path of
$\Pi_S^{(C)}( \cdot)$ corresponds to $P_C.$ More precisely, if
$P_C$ has an atom at time $S+t,$ say $\pi^{(C)}(S+t),$  we
coalesce $\Pi_S^{(C)}(t-)$ by $\pi^{(C)}(S+t)$ (using our new
labelling convention). For each $n$ we define
$$\text{dust}_n^{\text{coag}}(S,\cdot) := \text{dust} (\Pi_S^{(C)}(\cdot)) \cap
[n]$$ and $$ D_n^{\text{coag}}(S,\cdot) := \frac{1}{n} \#
\text{dust}_n^{\text{coag}}(S,\cdot).$$ We claim that for each $t
\ge 0$
$$\{i \le n : \forall s \in
[S,S+t] \; i \in \text{dust}_n(s) \} \subseteq
\text{dust}_n^{\text{coag}}(S,t).$$ Indeed suppose $j \in \{i \le
n : \forall s \in [S,S+t] \; i \in \text{dust}_n(s) \},$ then for
each $\pi^{(C)}(S+r)$ with $r \le t$ one has $j \in
\text{dust}(\pi^{(C)}(S+r))$ and hence $j$ has not yet coalesced
for the process $\Pi_S^{(C)}(\cdot).$ On the other hand, if there
exists a coalescence time $S+r$ such that $j \in
\text{dust}_n(S+r-)$ and $j \not \in \text{dust}_n(S+r)$ then it
is clear that $j$ also coalesces at time $S+r$ for
$\Pi^{(C)}_S(.)$ and hence $j \not \in
\text{dust}_n^{\text{coag}}(S,r).$ Thus we have that
$$D_n(S)-\frac{1}{n} \{i \le n : \forall s \in [S,S+t] \; i \in
\text{dust}_n(s) \} \le 1 - D_n^{\text{coag}}(S,t).$$ Now remark
that $$\{i \le n : \forall s \in [S,S+t] \; i \in \text{dust}_n(s)
\} \subseteq \text{dust}_n(S+t)$$ and thus
\begin{eqnarray*}
D_n(S) - D_n(S+t) &\le& D_n(S)-\frac{1}{n} \{i \le n : \forall s
\in [S,S+t] \; i \in \text{dust}_n(s) \} \\ &\le&
D_n(S)-\frac{1}{n} \{i \le n : \forall s \in [S,S+\theta] \; i \in
\text{dust}_n(s) \} \\ &\le& 1 - D_n^{\text{coag}}(S,\theta)
\end{eqnarray*}
(for the second inequality remark that $\{i \le n : \forall s \in
[S,S+t] \; i \in \text{dust}_n(s) \}$ is decreasing). We can now
apply the strong Markov property for the PPP $P_C$ at time $S$ and
we see that
\begin{eqnarray*}
P(1 - D_n^{\text{coag}}(S,\theta)>\epsilon) &=& P(1 -
D_n^{\text{coag}}(0,\theta)>\epsilon)
\\ &=& P\left( - \log
\left(D_n^{\text{coag}}(0,\theta) \right) > -\log (1-
\epsilon)\right) .
\end{eqnarray*}

Define $$\xi_n(t):=-\log (D_n^{\text{coag}}(0,t) ).$$ We know that
almost surely, for all $t \ge 0$ one has $\xi_n(t) \rightarrow
\xi(t)$ where $\xi(t)$ is a subordinator whose Lévy measure is
given by the image of $\nuc$ by the map $x \mapsto -\log (1-\sum_i
x_i).$ Hence, $P(\xi_n(\theta)> - \log(1-\epsilon)) \rightarrow
P(\xi(\theta)> - \log(1-\epsilon))$ when $n \rightarrow \infty.$
Thus, for any $\eta >0$ there exists a $\theta_1$ such that is
$\theta <\theta_1$ one has $\limsup_n P(\xi_n(\theta)> -
\log(1-\epsilon) <\eta.$ This bound being uniform in $S,$ the
conditions for applying Aldous' criterion are fulfilled.
\end{proof}

It is not hard to see that $D_n(\cdot),$ which takes its values in
$\{0,1/n,2/n,...,n/n\},$ is a Markov chain with the following
transition rates:
\begin{itemize}
\item if $k <n$ it jumps from $k/n$ to $(k+1)/n$ with rate $\ce n
(1-k/n)$,

\item if $k >0$ it jumps from $k/n$ to $r/n$ for any $r$ in
$0,....,k$ with rate $C_k^r \int_0^1 x^r (1-x)^{k-r}
\tilde{\nu}(dx)$ where $\tilde{\nu}$ is the image of $\nuc$ by the
map $\sfleche \mapsto [0,1] : x \mapsto (1-\sum_i x_i).$

\end{itemize}
Hence, if $A_n$ is the generator of the semi-group of $D_n$ one
necessarily has for any $f$ continuous
\begin{eqnarray} \label{genrateur An}
A_n f(k/n) &=& \frac{f((k+1)/n) - f(k/n)}{1/n}\ce(1-k/n) \\
\nonumber && + \sum_{r=1}^{k} (f((k-r)/n)  - f(k/n))C_k^r \int_0^1
x^r (1-x)^{k-r} \tilde{\nu}(dx).
\end{eqnarray}

We wish to define the $A_n$ so they will have a common domain,
hence we will let $A_n$ be the set of pairs of functions $f,g$
such that $f : [0,1] \mapsto \R$ is continuously differentiable on
$[0,1]$ and $f(D_n(t)) - \int_0^t g(D_n(s))ds$ is a martingale.
Note that continuously differentiable functions on $[0,1]$ are
dense in $C([0,1])$ the space of continuous functions on $[0,1]$
for the $L_{\infty}$ norm.

Hence $A_n$ is multivalued because for each function $f$, any
function $g$ such that $g(k/n)$ is given by (\ref{genrateur
An})will work. In the following we focus on the only such $g_n$
which is linear on each $[k/n,(k+1)/n].$

We know that $D_n \Rightarrow D$ in the space of càdlàg functions
and that $D_n$ is solution of the martingale problem associated to
$A_n$. Define $$Af(x):=f'(x)(1-x)\ce + \int_0^1 (f(\theta x)
-f(x)) \tilde \nu(d\theta).$$

In the following $\| f \| = \sup_{x \in [0,1]} |f(x)|.$
\begin{lemma}
One has $$\lim_{n \rightarrow \infty} \| g_n -Af \|=0.$$
\end{lemma}
\begin{proof}

We decompose $g_n$ into $g_n=g_n^{(1)} + g_n^{(2)}$ where both
$g_n^{(1)}$ and $g_n^{(2)}$ are linear on each $[k/n,(k+1)/n]$ and
$$g_n^{(1)}(k/n) = \frac{f((k+1)/n) -
f(k/n)}{1/n} \ce (1-k/n)$$ while $$g_n^{(2)}(k/n) = \sum_{r=1}^{k}
(f((k-r)/n)  - f(k/n))C_k^r \int_0^1 \theta^r (1-\theta)^{k-r}
\tilde{\nu}(d\theta).$$

One has that $\frac{f((k+1)/n) - f(k/n)}{1/n} \rightarrow f'(x)$
when $n\rightarrow \infty$ and $k/n \rightarrow x.$ Hence, as $f'$
is continuous on $[0,1],$ one has that$$\| g^{(1)}_n(x) - f'(x)
\ce(1-x) \| \rightarrow 0.$$

Let us now turn to the convergence of $g_n{(2)}.$ For a fixed $x$
and a fixed $\theta$ one has that
\begin{eqnarray*}
\sum_{r=1}^{[nx]} (f(r/n)  - f([nx]/n))C_{[nx]}^r \theta^r
(1-\theta)^{[nx]-r}  \rightarrow f(\theta x) - f(x)
\end{eqnarray*}
when $n \rightarrow \infty$ because $C_{[nx]}^r \theta^r
(1-\theta)^{[nx]-r} = P(B_{[nx],\theta}=r)$  where
$B_{[nx],\theta}$ is a $[nx],\theta$-binomial variable. We need
this convergence to be uniform in $x$. We proceed in two steps:
first it is clear that
$$\lim_n \sup_x (f(x) - f([nx]/n)) =0.$$

For the second part fix $\epsilon >0.$ There exists $\eta >0$ such
that $\forall x,y \in [0,1]$ one has $|x-y|\le \eta \Rightarrow
|f(x)-f(y)|<\epsilon.$

Next it is clear that there is a $n_0 \in \N$ such that $\forall
n\ge n_0$ and $\forall x \in [\eta,1]$ one has
\begin{eqnarray*}
&& P(B_{[nx],\theta} \in  [[nx](\theta -\eta),[nx](\theta +\eta)])
\\   \ge && P(B_{[n\eta],\theta} \in  [[n\eta](\theta
-\eta),[n\eta](\theta +\eta)] ) \\   > && 1-\epsilon.
\end{eqnarray*}

Hence, for $n \ge n_0$ and $x > \theta$
\begin{eqnarray*}
\sum_{r=1}^{[nx]} && f(r/n)C_{[nx]}^r \theta^r (1-\theta)^{[nx]-r}
\\&& \ge (1-\epsilon) \inf_{r \in [[nx](\theta -\eta),[nx](\theta
+\eta)] } f(r/n) \\ &&\ge (1-\epsilon) \inf_{\theta' \in [\theta -
\eta , \theta + \eta] } f(\frac{[nx]}{n}\theta')
\end{eqnarray*}
and
\begin{eqnarray*}
\sum_{r=1}^{[nx]} && f(r/n)C_{[nx]}^r \theta^r (1-\theta)^{[nx]-r}
\\ &&\le \sup_{r \in [[nx](\theta
-\eta),[nx](\theta +\eta)] } f(r/n) + \epsilon \|f\| \\
&& \le  \sup_{\theta'' \in [\theta - \eta , \theta + \eta] }
f(\frac{[nx]}{n}\theta'') + \epsilon \|f\|.
\end{eqnarray*}

Hence, for any $\epsilon'>0,$ by choosing $\epsilon$ and $\eta$
small enough, one can ensure that  there exists a $n_1$ such that
for all $n \ge n_1$ one has $$\sup_{x \in [\eta,1]} \left|
f(\theta x) - \sum_{r=1}^{[nx]}  f(r/n)C_{[nx]}^r \theta^r
(1-\theta)^{[nx]-r} \right| \le \epsilon'.$$

For $x <  \eta$ remark that
\begin{eqnarray*}
&& \left| f(\theta x) - \sum_{r=1}^{[nx]}  f(r/n)C_{[nx]}^r
\theta^r (1-\theta)^{[nx]-r} \right| \\ &\le & \left| f(\theta x)
- f(0)\right|+ \left| f(0) - \sum_{r=1}^{[nx]} f(r/n)C_{[nx]}^r
\theta^r (1-\theta)^{[nx]-r} \right|.
\end{eqnarray*}
We can bound $\sum_{r=1}^{[nx]} f(r/n)C_{[nx]}^r \theta^r
(1-\theta)^{[nx]-r}$ as follows:
\begin{eqnarray*}
&& P(B_{[nx],\theta} < [n\eta]+1 ) \inf_{s \le \eta}f(s) \\ &\le&
\sum_{r=1}^{[nx]} f(r/n)C_{[nx]}^r \theta^r (1-\theta)^{[nx]-r}
\\ &\le& P(B_{[nx],\theta} < [n\eta]+1 )\sup_{s \le \eta}f(s) + \|f\|
P(B_{[nx],\theta} \ge [n\eta]+1 ).
\end{eqnarray*}
Hence one has that
\begin{eqnarray*}
\lim_n \sup_x \left(\sum_{r=1}^{[nx]} (f(r/n))C_{[nx]}^r \theta^r
(1-\theta)^{[nx]-r} - f(\theta x) \right) =0.
\end{eqnarray*}
Finally we conclude that
\begin{eqnarray*}
\sup_x \left| \left[\sum_{r=1}^{[nx]} (f(r/n)  -
f([nx]/n))C_{[nx]}^r \theta^r (1-\theta)^{[nx]-r}\right] -\Bigg[
f(\theta x) - f(x) \Bigg] \right| \rightarrow 0
\end{eqnarray*}

We can then apply the dominated convergence theorem and we get
\begin{eqnarray*}
 \sup_x \left|g_n^{(2)}(x) - \int_0^1\left[ f(\theta x) - f(x) \right]
\tilde \nu (d\theta) \right| \rightarrow 0.
\end{eqnarray*}
Hence, one has $\|g_n^{(2)} - g^{(2)}\| \rightarrow 0$ where
$g^{(2)}(x)=\int_0^1 f(\theta x) - f(x) \tilde{\nu}(d\theta).$

\end{proof}
One can now use Lemma 5.1 in \cite{ethier_kurtz} to see that $D$
must be solution of the Martingale Problem associated to $A$.
Hence one can use Theorem III.2.26 in \cite{jacod_shiryaev} to see
that $D$ is solution of $$dD(t) = \ce (1-D(t))dt + \int_0^1
D(t-)(\theta-1)p(dt,d\theta)$$ where $p(dt,d\theta)$ is the
counting measure for the PPP $(t,1-\sum_i x_i)$ with measure
intensity $\tilde \nu.$

Recall that $\xi(t)=-\log (D(t)).$ By taking $f = g \circ -\log$
where $g$ is such that $f \in \mathcal D(A),$ and using standard
results (see \cite{jacod_prob_de_mg}) one has that $\xi$ is
solution of the martingale problem associated to the generator
$$A'g(x) = -\ce (e^x-1) g'(x) + \int_0^1 (g( x - \log \theta)
-g(x)) \tilde \nu(d\theta).$$

Hence it is easily seen that $\xi$ is solution of the SDE
$$d\xi(t) = d\sigma(t) - \ce (e^{\xi(t)}-1)dt$$ where $\sigma$ is
a drift-free subordinator whose Lévy measure is the image of
$\tilde \nu$ by $x \mapsto -\log x.$

\section{Proofs}

\subsection{Proof of Proposition \ref{prop taux de sauts}}

The compatibility of the chains $\Pi_{|[n]}$ can be expressed in
terms of transition rates as follows: For $m <n \in \N$  and
$\pi,\pi' \in \P_{n}$ one has
$$q_m(\pi_{|[m]},\pi'_{|[m]})=\sum_{\pi'' \in\P_{n} :
\pi''_{|[m]}=\pi'_{|[m]}} q_{n}(\pi,\pi'').$$

Consider $\pi \in \P_n$ such that $\pi=(B_1,B_2,...,B_m,\o,...)$
has $m \le n$ non-empty blocks. Call $w_i = \inf\{k \in B_i\}$ the
least element of $B_i$ and $\sigma$ a permutation of $[n]$ that
maps every $i \le m $ on $w_i$. Let $\pi'$ be an element of
$\P_m$, then the restriction of the partition
$\sigma(\coag(\pi,\pi'))$ to $[m]$ is given by: for $i,j \le m$
\begin{eqnarray*}
i \overset{\sigma(\coag(\pi,\pi'))}{\sim} j &\Leftrightarrow&
\sigma(i) \overset{\coag(\pi,\pi')}{\sim} \sigma(j) \\
&\Leftrightarrow& \exists k,l : \sigma(i) \in B_k, \sigma(j) \in
B_l , k\overset{\pi'}{\sim} l \\ &\Leftrightarrow& i
\overset{\pi'}{\sim} j
\end{eqnarray*}
and hence
\begin{eqnarray}\label{toto1}
\sigma(\coag(\pi,\pi'))_{|[m]}=\pi'.
\end{eqnarray}
By definition $C_n(\pi,\pi')$ is the rate at which the process
$\sigma(\Pi_{|[n]}(\cdot))$ jumps from $\sigma(\pi)$ to
$\sigma(\coag(\pi,\pi')).$  Hence, by exchangeability
$$C_n(\pi,\pi') = q_n(\sigma(\pi),\sigma(\coag(\pi,\pi'))).$$

Remark that $\sigma(\pi)_{|[m]}=\0_m.$ Hence if $\pi''$ is a
coalescence of $\sigma(\pi)$ it is completely determined by
$\pi''_{|[m]}.$ Thus, for all $\pi'' \in \P_n$ such that
$\pi_{|[m]}''=\sigma(\coag(\pi,\pi'))_{|[m]}$ and $\pi'' \neq
\sigma(\coag(\pi,\pi'))$  one has
\begin{eqnarray}\label{taux}
q_n(\sigma(\pi),\pi'') =0.
\end{eqnarray}
For each $\pi \in \P_n$ define $$Q_n(\pi,m):=\{\pi' \in \P_n :
\pi_{|[m]}=\pi'_{|[m]}\}$$ (for $\pi \in \P$ we will also need
$Q(\pi,m):=\{\pi' \in \P : \pi_{|[m]}=\pi'_{|[m]}\}$). Clearly,
(\ref{taux}) yields
$$q_n(\sigma(\pi),\sigma(\coag(\pi,\pi')) = \sum_{\pi''
\in Q_n(\sigma(\coag(\pi,\pi')),m)}q_n(\sigma(\pi),\pi'')$$
because there is only one non-zero term in the right hand-side
sum. Finally recall (\ref{toto1}) and use the compatibility
relation to have
\begin{eqnarray*}
C_n(\pi,\pi')&=&q_n(\sigma(\pi),\sigma(\coag(\pi,\pi')) \\ &=&
\sum_{\pi'' \in Q_n(\sigma(\coag(\pi,\pi')),m)}
q_n(\sigma(\pi),\pi'')
\\&=& q_m(\sigma(\pi)_{|[m]},\sigma(\coag(\pi,\pi'))_{|[m]}) \\&=&
q_m(\0_m,\pi') \\&=& C_m(\0_m,\pi') \\ &:=& C_m(\pi').
\end{eqnarray*}

Let us now take care of the fragmentation rates. The argument is
essentially the same as above. Suppose $B_k =
\{n_1,...,n_{|B_k|}\}$. Let $\sigma$ be a permutation of $[n]$
such that of all $j \le |B_k| $ one has $\sigma(j)=n_j$. Hence, in
$\sigma(\pi)$ the first block is $[|B_k|]$. The process
$\sigma(\Pi_{|[n]}(\cdot))$ jumps from $\sigma(\pi)$ to the state
$\sigma(\frag(\pi,\pi',k))$ with rate $F_n(\pi,\pi',k)$. Remark
that for $i,j \le |B_k|$
\begin{eqnarray*}
i\overset{\sigma(\frag(\pi,\pi',k))}{\sim}j &\Leftrightarrow&
\sigma(i) \overset{\frag(\pi,\pi',k)}{\sim}\sigma(j) \\
&\Leftrightarrow& n_i \overset{\frag(\pi,\pi',k)}{\sim}n_j \\
&\Leftrightarrow& i \overset{\sigma(\pi')}{\sim} j
\end{eqnarray*}
and hence
\begin{eqnarray}\label{toto2}
\sigma(\frag(\pi,\pi',k)) =  \frag(\sigma(\pi),\sigma(\pi'),1).
\end{eqnarray}
Thus by exchangeability
$F_n(\pi,\pi',k)=F_n(\sigma(\pi),\sigma(\pi'),1)$, and it is
straightforward to see that by compatibility that
$$F_n(\sigma(\pi),\sigma(\pi'),1)=F_{|B_k|}(\1_{|B_k|},\sigma(\pi'),1) = F_{|B_k|}(\sigma(\pi')).$$

The invariance of the rates $C_n(\0_n,\pi')$ and
$F_n(\1_n,\pi',1)$ by permutations of $\pi'$ is also a direct
consequence of exchangeability. In particular
$F_{|B_k|}(\sigma(\pi'))=F_{|B_k|}(\pi')$ and thus we conclude
that $F_n(\pi,\pi',k)=F_{|B_k|}(\pi')$.

\subsection{Proof of Theorem \ref{Markov pour X}}

We first have to introduce a few notations: let $B(i,t)$ denote
the block that contains $i$ at time $t$ and define
\begin{itemize}
\item $\bar \lambda_i(t)=\bar \lambda_{B_i(t)}$ and $\underline
\lambda_i(t)=\underline \lambda_{B_i(t)},$ \item $\bar
\lambda(i,t)=\bar \lambda_{B(i,t)}$ and $\underline
\lambda(i,t)=\underline \lambda_{B(i,t)}.$
\end{itemize}

In the following we will use repeatedly a coupling technique that
can be described as follows: Suppose $\Pi$ is an EFC process
constructed with the PPP $P_F$ and $P_C$, we choose $T$ a stopping
time for $\Pi,$ at time $T$ we create a fragmentation process
$(\Pi^{(F)}(T+s),s\ge 0)$ started from $\Pi^{(F)}(T)=\Pi(T)$ and
constructed with the PPP $(P_F(T+s),s\ge 0)$. We call
$(B_1^{(F)}(T+s),B_2^{(F)}(T+s),...)$ the blocks of
$\Pi^{(F)}(T+s)$ and $\bar \lambda_i^{(F)}(T+s), \underline
\lambda_i^{(F)}(T+s)$ the corresponding limsup and liminf for the
frequencies. The processes $\Pi(T+s)$ and $\Pi^{(F)}(T+s)$ are
coupled. More precisely, remark that for instance
$$B^{(F)}_1(T+s) \subseteq B_1(T+s) ,  \forall s\ge 0$$ because if
$i \in B^{(F)}_1(T+s)$ it means that there is no $r \in [T,T+s]$
such that $k(r)=1$ and $1 \overset{\pi^{(F)}(r)}{ \not \sim} i$
and hence $i\in B_1(T+s).$

For any exchangeable variable $\Pi=(B_1,B_2,...)$, and $A \subset
\N$ independent of $\Pi$ one can easily see that almost surely for
each $i \in \N$
$$\limsup_{n \rightarrow \infty} \frac{\#\{ k
\le n : k \in B_i \cap A \}}{n}= \bar \lambda_A \|B_i\|$$ and
$$\liminf_{n \rightarrow \infty} \frac{\#\{ k \le n : k \in B_i
\cap A \}}{n}= \underline \lambda_A \|B_i\|.$$ Hence, if we start
a homogeneous fragmentation $(\Pi^{(F)}(T+s),s\ge 0)$ from a
partition that does not necessarily admit asymptotic frequencies,
say $\Pi^{(F)}(T)=(...,A,....)$ (i.e., $A$ is one of the block in
$\Pi^{(F)}(0)$), we still have that if $a$ designates the least
element of $A$ then almost surely
\begin{eqnarray}
\label{conv des limsups} \bar \lambda^{(F)}(a,T+s) \rightarrow
\bar \lambda_A
\end{eqnarray} and
$$\underline \lambda^{(F)}(a,T+s) \rightarrow \underline \lambda_A$$
when $s \searrow 0.$

To prove Theorem \ref{Markov pour X}, it suffices to prove the
existence of the asymptotic frequency of $B_1(t)$ simultaneously
for all $t$, the same demonstration then apply to the $B(i,t)$ for
each $i$. As $\Pi(t)$ is an exchangeable process we already know
that $\|B_1(q)\|$ exists simultaneously for all $q \in \Q.$ For
such $q$ we thus have that $\bar\lambda_1(q)= \underline
\lambda_1(q).$ Hence, it suffices to show that $\bar\lambda_1(t)$
and $\underline \lambda_1(t)$ are both càdlàg processes. In the
following we write $q \searrow \! \! \searrow t$ or $q \nearrow \!
\! \nearrow t$ to mean $q$ converges to $t$ in $\Q$ from below
(resp. from above).

The first step is to show that:
\begin{lemma}
Almost surely, the process  $(L(t),t \ge0)$ defined by
\begin{eqnarray*}
\forall t \ge 0 : L(t) := \lim_{q \searrow \! \! \searrow t} \bar
\lambda_1(q) = \lim_{q \searrow \! \! \searrow t} \underline
\lambda_1(q)
\end{eqnarray*}
exists and is càdlàg.
\end{lemma}

\begin{proof}
Using standard results (see for instance \cite[Theorem
62.13]{rogers_william}), and recalling that $\bar \lambda_1$ and
$\underline \lambda_1$ coincide on $\Q$, one only need to show
that $q \mapsto \bar \lambda_1(q)=\underline \lambda_1(q)$ is a
regularisable process, that is
\begin{eqnarray*}
\lim_{q \searrow  \! \! \searrow t} \bar \lambda_1(q) = \lim_{q
\searrow  \! \! \searrow t} \underline \lambda_1(q) \text{ exist
for every real } t\ge 0,
\\\lim_{q \nearrow \! \! \nearrow t} \bar \lambda_1(q) = \lim_{q
\nearrow \! \! \nearrow t} \underline \lambda_1(q) \text{ exist
for every real } t\ge 0. \end{eqnarray*}

Using \cite[Theorem 62.7]{rogers_william}, one only has to verify
that whenever $N\in\N$ and $a,b \in \Q$ with $a<b$, almost surely
we have
\begin{eqnarray*}
\sup \{ \bar \lambda_1(q) : q \in \Q^+ \cap [0,N] \} = \sup \{
\underline \lambda_1(q) : q \in \Q^+ \cap [0,N] \}<\infty
\end{eqnarray*}
and
\begin{eqnarray*}
U_N(\bar \lambda_1;[a,b]) = U_N(\underline \lambda_1;[a,b])<\infty
\end{eqnarray*}
where $U_N(\bar \lambda_1;[a,b])$ is the number of upcrossings of
$\bar \lambda_1$ from $a$ to $b$ during $[0,N].$ By definition
$\sup \{ \bar \lambda_1(q) : q \in \Q^+ \cap [0,N] \} \le 1$ and
$\sup \{ \underline \lambda_1(q) : q \in \Q^+ \cap [0,N] \} \le
1.$ Suppose that $q \in \Q$ is such that $\bar \lambda_1(q)>b.$
Then if we define $s=\inf \{r \ge 0 : \bar \lambda_1(q+r) \le a\}$
one can use the Markov property and the coupling with a
fragmentation $(\Pi^{(F)}(q+r),r \ge 0)$ started from $\Pi(q)$,
constructed with the PPP $(P_F(q+r),r\ge 0)$ to see that $s\ge
\theta$ where $\theta$ is given by $\theta := \inf \{r \ge 0 :
\bar \lambda^{(F)}_1(t+r) \le a\}.$ If one has a sequence
$L_1<R_1<L_2<R_2,....$ in $\Q$ such that $\bar
\lambda_1(L_i)<a<b<\bar \lambda_1(R_i),$ then one has that for
each $i,$ $R_i-L_i>\theta_i$ where $(\theta_i)_{i \in \N}$ is an
i.i.d. sequence with same distribution as $\theta$. Hence
$P(U_N(\bar \lambda_1;[a,b])=\infty)=0.$
\end{proof}

The next step is the following:
\begin{lemma}
Let $T$ be a stopping time for $\Pi.$ Then one has $$\sum_{i \in
\N} \bar \lambda_i(T) \le 1$$ and $\bar{\lambda}_1$ and
$\underline \lambda_1$ are right continuous at $T$.
\end{lemma}
\begin{proof}

For the first point, suppose that $\sum_i \bar \lambda_i(T) =1+
\gamma>1.$ Then there exists $n\in \N$ such that $\sum_{i \le N}
\bar \lambda_i(T) >1+\gamma/2$. Call $w_i(t)$ the least element of
$B_i(t).$ Let $S$ be the random stopping time defined as the first
time after $T$ such that $S$ is a coalescence involving at least
two of the $w_N(T)$ first blocks $$S = \inf \{s \ge T :
\pi^{(C)}(s) \not \in Q(\0,w_N(T)) \}.$$ Hence, between $T$ and
$S,$ for each $i \le N$ one has that $w_i(T)$ is the least element
of its block. Applying the Markov property in $T$ we have that
$S-T$ is exponential with a finite parameter and is thus almost
surely positive.

Define $(\Pi^{(F)}(T+s),s\ge 0)$ as the fragmentation process
started from $\Pi(T)$ and constructed from the PPP $(P_F(T+s),s
\ge 0).$ On the time interval $[T,S]$ one has that for each $i$,
the block of $\Pi^{(F)}$ that contains $w_i$ is included in the
block of $\Pi$ that contains $w_i$ (because the last might have
coalesced with blocks whose least element is larger than
$w_N(T)$).

Fix $\epsilon >0,$ using (\ref{conv des limsups}) and the above
remark, one has that for each $i \le N$ there exists a
$\theta_i>0$ such that for all $t \in [T,T+\theta_i]$ one has
$$\bar \lambda(w_i(T),t)>(1-\epsilon)\bar \lambda(w_i(T),T).$$ Thus,
if $\theta = \min \theta_i$ one has that $$\min_{s \in
[T,T+\theta]} \sum_i \bar \lambda_i(s)
>(1+\gamma/2)(1-\epsilon).$$ Choosing $\epsilon$ small enough yields a
contradiction with the fact that almost surely for all $t \in \Q$
one has $\sum \bar \lambda_i(t)\le 1.$

Fix $\epsilon>0,$ the first part of the lemma implies that there
exists $N_{\epsilon} \in \N$ such that
$$\sum_{i \ge N_{\epsilon}} \bar \lambda_i(T) \le \epsilon.$$ Let
$(\Pi^{(F)}(T+s),s\ge0)$ be as above a fragmentation started from
$\Pi(T)$ and constructed with the PPP $(P_F(T+s), s\ge 0).$ As we
have remarked
\begin{eqnarray}
\label{borne inf pour fa sup}\bar \lambda_1(T+s) \ge \bar
\lambda_1^{(F)}(T+s) \rightarrow  \bar \lambda_1^{(F)}(T) \\
\label{borne inf pour fa inf} \underline \lambda_1(T+s) \ge
\underline \lambda_1^{(F)}(T+s)\rightarrow  \underline
\lambda_1^{(F)}(T)
\end{eqnarray}
when $s \searrow 0.$

Now consider $$S=\inf \{s \ge T : \pi^{(C)}_{|[N_{\epsilon}]}(s)
\neq \0_{N_{\epsilon}} \}$$ the first coalescence time after $T$
such that $\pi^{(C)}_{|[N_{\epsilon}]}(s) \neq \0_{N_{\epsilon}}.$
One has $\forall s \in [T,S]$
\begin{eqnarray*}
\bar \lambda_1^{(F)}(s) &\le&  \bar \lambda_1^{(F)}(T) + \sum_{i
\ge N_{\epsilon}} \bar \lambda_i(T) \le \bar \lambda_1^{(F)}(T) +
\epsilon
\\ \underline \lambda_1^{(F)}(T+s) &\le& \underline
\lambda_1^{(F)}(T)+ \sum_{i \ge N_{\epsilon}} \bar \lambda_i(T)
\le \underline \lambda_1^{(F)}(T) + \epsilon.
\end{eqnarray*}
Thus $\bar{\lambda}_1(T+s) \rightarrow \bar \lambda_1(T)$ and
$\underline \lambda_1(T+s) \rightarrow \underline \lambda_1(T)$
when $s\searrow 0.$
\end{proof}

To conclude the demonstration of the first point of Theorem
\ref{Markov pour X}, remark that as the map $\Pi(t)\mapsto \bar
\lambda_1(t)$ is measurable in $\F_t,$ the right-continuous usual
augmentation of the filtration, one has that for any $\epsilon >0$
$$\inf \{ t : |\limsup_{s \searrow 0} \bar \lambda_1(t+s)-\bar
\lambda_1(t)|>\epsilon \}$$ or $$\inf \{ t : |\liminf_{s \searrow
0} \bar \lambda_1(t+s)-\bar \lambda_1(t)|>\epsilon \}$$ are
stopping times for $\Pi$ in $\F.$ The above lemma applies and
hence this stopping times are almost surely infinite. The same
argument works for $\underline \lambda_1$. This shows that $\bar
\lambda_1$ and $\underline \lambda_1$ are almost surely
right-continuous processes. As they coincide almost surely with
$L$ on the set of rationals, they coincide everywhere and hence
their paths are almost surely càdlàg.

Before we can prove rigourously that $X(t)$ is a Feller process,
as stated in Theorem \ref{Markov pour X}, we have to pause for a
moment to define a few notions related to the laws of EFC
processes conditioned on their starting point. By our definition,
an EFC process $\Pi$ is exchangeable. Nevertheless, if $P$ is the
law of $\Pi$ and $P_{\pi}$ is the law of $\Pi$ conditionally on
$\Pi(0)=\pi,$ one has that as soon as $\pi\neq \0$ or $\1,$ the
process $\Pi$ is not exchangeable under $P_{\pi}$ (because for
instance $\Pi(0)$ is not exchangeable). The process $\Pi$
conditioned by $\Pi(0)=\pi$ (i.e., under the law $P_{\pi}$) is
called an EFC evolution. Clearly one can construct every EFC
evolution exactly as the EFC processes, or more precisely, given
the PPP's $P_F$ and $P_C$ one can then choose any initial state
$\pi$ and construct the EFC evolution $\Pi, \Pi(0)=\pi$ with $P_F$
and $P_C$ as usually. Let us first check quickly that under
$P_{\pi}$ we still have the existence of $X(t)$ simultaneously for
all $t.$

In the following we will say that a partition $\pi \in \P$ is
\textit{good} if $\Lambda(\pi)$ exists, there are no finite blocks
of cardinal greater than 1 and either $\text{dust}(\pi)=\o$ or
$\|\text{dust}(\pi)\| >0.$

\begin{lemma}
For each $\pi \in \P$ such that $\pi$ is good, then $P_{\pi}$-a.s.
the process $X(t)=\Lambda(\Pi(t))$ exists for all $t$
simultaneously and we call $Q_{\pi}$ its law.
\end{lemma}

\begin{proof}
Consider $\pi=(B_1,B_2,...)$ a good partition. For each $i \in \N$
such that $\# B_i =\infty,$ let $f_i : \N \mapsto \N$ be the only
increasing map that send $B_i$ on $\N.$ Let $B_0 = \cup_{i : \#
B_i <\infty} B_i$ and if $B_0$ is infinite(which is the case
whenever it is not empty) set $g : \N \mapsto \N$ the unique
increasing map that send $B_0$ onto $\N.$

Using the exchangeability properties attached to the PPP's $P_F$
and $P_C$ one can easily see that for each $i \in \N$ such that
$\# B_i =\infty,$ $$ f_i(\Pi(t) \cap B_i)   $$ and $$ g(\Pi(t)
\cap B_0)$$ are EFC processes with initial state $\1$ for the
first ones and $\0$ for the later. Hence for each $i$ one has that
$ f_i(\Pi(t) \cap B_i) $ has asymptotic frequencies
$X^{(i)}(t):=\Lambda(f_i(\Pi(t) \cap B_i))$ simultaneously for all
$t.$ Thus it is not hard to see from this that $\Pi(t) \cap B_i$
has asymptotic frequencies simultaneously for all $t,$ namely
$\|B_i\| X^{(i)}(t).$

Fix $\epsilon >0,$ there exists $N_{\epsilon}$ such that
$$\|B_0 \| +\sum_{i \le \N_{\epsilon}} \|B_i\| \ge
1-\epsilon.$$ If we call $\Pi(t)=(B_1(t),B_2(t),...)$ the blocks
of $\Pi(t),$ we thus have that for $j\in \N$ fixed $$\bar
\lambda_j(t) \le  \sum_{i \le N_{\epsilon}} \|B_j(t) \cap B_i\| +
\epsilon$$ and $$\underline \lambda_j(t) \ge  \sum_{i \le
N_{\epsilon}} \|B_j(t) \cap B_i\|.$$ Hence
$$\sup_{t \ge 0} \sup_{i \in \N}
(\bar \lambda_i(t)-\underline \lambda_i(t) )\le \epsilon.$$ As
$\epsilon$ is arbitrary this shows that almost surely $\sup_{t \ge
0} \sup_{i \in \N} (\bar \lambda_i(t)-\underline \lambda_i(t)) =
0.$ We call $Q_{\pi}$ the law of $X(t)$ under $P_{\pi}.$
\end{proof}

Although EFC evolutions are not exchangeable, they do have a very
similar property:

\begin{lemma}\label{exchangeability for evolution}
Let $(\Pi_1(t), t\ge 0)$ be an EFC evolution with law $P_{\pi_1}$
(i.e., $P(\Pi_1(0)=\pi_1)=1$) and with characteristics
$\nuf,\nuc,\cc$ and $\ce.$ Then for any bijective map $\sigma : \N
\mapsto \N$ the process $\Pi_2(t):=(\sigma^{-1}(\Pi_1(t)),t \ge
0)$ is an EFC evolution with law $P_{\sigma^{-1}(\pi_1)}$ and same
characteristics.
\end{lemma}
\begin{proof}

Consider $\Pi_1(t)=(B^{(1)}_1(t),B^{(1)}_2(t),...)$ an EFC
evolution with law $P_{\pi_1}$ (i.e., started from $\pi_1$) and
constructed with the PPP's $P_F$ and $P_C.$ Let $\pi_2 =
\sigma^{-1}(\pi_1)$ and $(\Pi_2(t),t\ge
0)=(\sigma^{-1}(\Pi_1(t)),t\ge 0).$ For each $t \ge 0$ and $k \in
\N$ call $\phi(t,k)$ the label of the block
$\sigma(B^{(1)}_k(t-))$ in $\Pi_2(t-).$ By construction,
$\Pi_2(t)$ is a $\P$-valued process started from $\pi_2.$ When
$P_F$ has an atom, say $(k(t),\pi^{(F)}(t))$ the block of
$\Pi_2(t-)$ which fragments has the label $\phi(t,k(t))$ and the
fragmentation is done by taking the intersection with
$\sigma^{-1}(\pi^{(F)}(t)).$ Call $\tilde P_F$ the point process
of the images of the atoms of $P_F$ by the transformation
$$(t,k(t),\pi^{(F)}(t)) \mapsto
(t,\phi(t,k(t)),\sigma^{-1}(\pi^{(F)}(t))).$$ If $t$ is an atom
time for $P_C$, say $\pi^{(C)}(t),$ then $\Pi_2$ also coalesces at
$t$, and if the blocks $i$ and $j$ merge at $t$ in $\Pi_1$ then
the blocks $\phi(t,i)$ and $\phi(t,j)$ merge at $t$ for $\Pi_2$,
hence the coalescence is made with the usual rule by the partition
$\phi^{-1}(t,\pi^{(C)}(t)).$ Call $\tilde P_C$ the point process
image of $P_C$ by the transformation
$$(t,\pi^{(C)}(t)) \mapsto (t,\phi^{-1}(t,\pi^{(C)}(t))).$$

We now show that $\tilde P_C$ and $\tilde P_F$ are PPP with the
same measure intensity as $P_C$ and $P_F$ respectively. The idea
is very close to the proof of Lemma 3.4 in \cite{moi_1}. Let us
begin with $\tilde P_F.$ Let $A \subset \P$ such that
$(\mu_{\nuf}+\ce \meros)(A)<\infty$ and define $$N_A^{(i)}(t):=\#
\{ u \le t : \sigma(\pi^{(F)}(u)) \in A, k(u)=i \}.$$ Then set
$$N_A(t) := \# \{u \le t :\sigma(\pi^{(F)}(u)) \in A,
\phi(u,k(u))=1 \}.$$ By definition $$d N_A(t)= \sum_{i=1}^{\infty}
\indic{\phi(t,i)=1}dN^{(i)}_A(t).$$ The process $N_A$ is
increasing, càdlàg and has jumps of size $1$ because by
construction the $N_A^{(i)}$ do not jump at the same time almost
surely. Define the counting processes $\bar N^{(i)}_A(t)$ by the
following differential equation $$ d \bar N^{(i)}_A(t) =
\indic{\phi(t,i)=1}dN^{(i)}_A(t).$$ It is clear that
$\indic{\phi(t,i)=1}$ is adapted and left-continuous in $(\F_t)$
the natural filtration of $\Pi_1$ and hence predictable. The
$N_A^{(i)}(\cdot)$ are i.i.d. Poisson process with intensity
$(\mu_{\nuf}+ \ce \meros)(A)=(\mu_{\nuf}+ \ce
\meros)(\sigma^{-1}(A))$ in $(\F_t).$ Thus for each $i$ the
process
\begin{eqnarray*}
 M_A^{(i)}(t)&=&\bar
N_A^{(i)}(t)-(\mu_{\nuf}+\ce \meros)(A)\int_0^t
\indic{\phi(u,i)=1} du
\\ &=& \int_0^t  \indic{\phi(u,i)=1} d\big(N_A^{(i)}(u)
-(\mu_{\nuf}+\ce \meros)(A)u\big)
\end{eqnarray*}
is a square-integrable martingale.

Define $$M_A(t):=\sum_{i=1}^{\infty} \int_0^t  \indic{\phi(u,i)=1}
d(N_A^{(i)}(u) -(\mu_{\nuf}+\ce \meros)(A)u).$$ For all $i\neq j$
, for all $t\ge 0$ one has
$\indic{\phi(t,i)=1}\indic{\phi(t,j)=1}=0$ and for all $t \ge 0$
one has $\sum_{i=1}^{\infty} \indic{\phi(u,i)=1}=1,$ the
$M_A^{(i)}$ are orthogonal (because they do not share any
jump-time) and hence the oblique bracket of $M_A$ is given by
\begin{eqnarray*}
<M_A>(t) &=& \sum_{i=1}^{\infty} \left\langle \int_0^t
\indic{\phi(u,i)=1} d(N_A^{(i)}(u) -(\mu_{\nuf}+\ce \meros)(A)u)
\right\rangle \\ &=& \mu_{\nuf}(A)t.
\end{eqnarray*}
Hence $M_A$ is a $L_2$ martingale. This shows that $N_A(t)$ is
increasing càdlàg with jumps of size 1 and has $(\mu_{\nuf}+\ce
\meros)(A) t$ as compensator. We conclude that $N_A(t)$ is a
Poisson process of intensity $(\mu_{\nuf}+\ce \meros)(A).$ Now
take $B \subset \P$ such that $A\cap B =\o$ and consider $N_A(t)$
and $N_B(t)$, clearly they do not share any jump time because the
$N_A^{(i)}(t)$ and $N_B^{(i)}(t)$ don't. Hence
$$P_F^{(1)}(t)=\{\sigma(\pi^{(F)}(u)): u\le t , \phi(u,k(t))=1
\}$$ is a PPP with measure-intensity $(\mu_{\nuf} +\ce \meros ).$
Now, by the same arguments $$P_F^{(2)}(t)=\{\sigma(\pi^{(F)}(u)):
u\le t , \phi(u,k(t))=2 \}$$ is also a PPP with measure-intensity
$(\mu_{\nuf}+\ce\meros)$ independent of $P_F^{(1)}.$ By iteration
we see that $\tilde P_F$ is a PPP with measure intensity
$(\mu_{\nuf}+\ce \meros) \otimes \#.$

Let us now treat the case of $\tilde P_C.$ The main idea is very
similar since the first step is to show that for $n\in\N$ fixed
and $\pi \in \P$ such that $\pi_{|[n]} \neq \0_n$ one has that the
counting process $$N_{\pi,n}(t) = \#\{ u \le t :
\phi^{-1}(u,\pi^{(C)}(u))_{|[n]} = \pi_{|[n]} \}$$ is a Poisson
process with intensity $(\mu_{\nuc}+ \cc \kappa)(Q(\pi,n)).$

For each unordered collection of $n$ distinct elements in $\N$,
say $\bold a = a_1,a_2,...,a_n,$  let $\sigma_{\bold a}$ be a
permutation such that for each $i \le n,$ $\sigma_{\bold
a}(i)=a_i.$

For each $\bold a$ define $$N_{\bold a,\pi}(t) = \#\{u \le t :
(\sigma_{\bold a}(\pi^{(C)}(u)))_{|[n]}=\pi_{|[n]} \}.$$ By
exchangeability $N_{\bold a,\pi}(t)$ is a Poisson process with
measure intensity $(\mu_{\nuc}$ $+ \cc \kappa)(Q(\pi,n)).$

By construction $$d N_{\pi,n}(t) = \sum_{\bold a} \prod_{i=1}^n
\indic{\phi(t,a_i)=i} d N_{\bold a,\pi}(t).$$ We see that we are
in a very similar situation as before: the $N_{\bold a,\pi}(t)$
are not independent but at all time $t$ there is exactly one
$\bold a$ such that $\prod_{i=1}^n \indic{\phi(t,a_i)=i} =1$ and
hence one can define orthogonal martingales $M_{\bold a}(t)$ as we
did for the $M_A^{(i)}(t)$ above and conclude in the same way that
$N_{\pi,n}(t)$ is a Poisson process with measure intensity
$(\mu_{\nuc}+ \cc \kappa)(Q(\pi,n)).$ If we now take $\pi' \in \P$
such that $\pi'_{|[n]} \neq \pi_{|[n]}$ we have that
$N_{\pi,n}(t)$ and $N_{\pi',n}(t)$ are independent because for
each fixed $\bold a$ the processes given by the equation $$d
M_{\bold a , \pi}(t) =\prod_{i=1}^n \indic{\phi(t,a_i)=i} d
N_{\bold a,\pi}(t)$$ and
$$d M_{\bold a , \pi'}(t) =\prod_{i=1}^n \indic{\phi(t,a_i)=i} d
N_{\bold a,\pi'}(t)$$ respectively does not have any common jumps.
Hence $N_{\pi,n}(t)$ and $N_{\pi',n}(t)$ are independent and thus
we conclude that $\tilde P_C$ is a PPP with measure intensity
$\mu_{\nuc} + \cc \kappa.$

Putting the pieces back together we see that $\Pi_2$ is an EFC
evolution with law $P_{\pi_2}$ and same characteristics as
$\Pi_1.$

\end{proof}

For each $\pi \in \P$ such that $\Lambda(\pi)=x$ exists, and for
each $k \in \N$ we define $n_{\pi}(k)$ the label of the block of
$\pi$ which corresponds to $x_k$, i.e., $\| B_{n_{\pi}(k)}
\|=x_k.$ In the case where two $B_k$'s have the same asymptotic
frequency we use the order of the least element, i.e., if there is
$i$ such that $x_i =x_{i+1}$ one has $n_{\pi}(i)<n_{\pi}(i+1).$
The map $i \mapsto n_{\pi}(i)$ being bijective, call $m_{\pi}$ its
inverse. Furthermore we define $B_0 = \cup_{i :\|B_i\|=0} B_i$ and
$x_0 =\|B_0\| =  1 -\sum_{i \in \N} x_i.$ In the following we will
sometimes write $\pi=(B_0,B_1,...)$ and $x=(x_0,x_1,...).$

Let $\pi=(B_1,B_2,...), \pi' =(B'_1,B'_2,..) \in \P$ be two good
partitions. We write $\Lambda(\pi)=x=(x_1,x_2,...)$ and
$\Lambda(\pi') = x'=(x'_1,x'_2,...).$ Suppose furthermore that
either $x_0=0$ and $x'_0=0$ or they are both strictly positive and
that
$$\inf \{ k \in \N : x_k=0 \} = \inf \{ k \in \N : x'_k=0 \}.$$
Define $\sigma_{\pi,\pi'}$ the unique bijection $\N \mapsto \N$
that map every $B_{n_{\pi}(i)}$ onto $B'_{n_{\pi'}(i)}$ such that
if $i,j \in B_{n_{\pi}(i)}$ with $i<j$ then $\sigma(i)<\sigma(j).$
Note that this definition implies that
$\sigma_{\pi,\pi'}(B_0)=B'(0)$. Furthermore we have $\pi'
=\sigma_{\pi,\pi'}^{-1}(\pi).$

We will use the following technical lemma:
\begin{lemma}\label{technique 1}
For $\pi,\pi'$ fixed in $\P$ verifying the above set of
hypothesis, let $\Pi(t)=(B_1(t),B_2(t),...)$ be an EFC evolution
started from $\pi$ with law $P_{\pi}$, then
\begin{itemize}
\item $\Lambda(\pi \cap \Pi(t))$  exists almost surely for all $t$
simultaneously where $\pi \cap \Pi(t)$ is defined by $i
\overset{\pi \cap \Pi(t)}{\sim}j$ if and only if we have both
$i\overset{\Pi(t)}{\sim}j$ and $i\overset{\pi}{\sim}j.$

\item $\Lambda(\sigma_{\pi,\pi'}^{-1}(\Pi(t) \cap \pi))$ also
exists a.s. for all $t\ge 0$ and for each $j,k \in \N$ one has
$$\|\sigma_{\pi,\pi'}( B_j(t)\cap B_k) \| =
\frac{x'_{m_{\pi}(k)}}{x_{m_{\pi}(k)}} \| B_j(t)\cap B_k \|.$$

\item $\Lambda(\sigma_{\pi,\pi'}^{-1}(\Pi(t)))$ exists a.s.
$\forall t \ge 0$ and for each $j$
$$\|\sigma_{\pi,\pi'}(B_j(t)) \| =\sum_{k \ge 0} \frac{x'_{m_{\pi}(k)}}{x_{m_{\pi}(k)}} \|
B_j(t)\cap B_k  \|.$$
\end{itemize}
\end{lemma}
\begin{proof}
For $B \subset \N$ call $F_B$ the increasing map that send $\N$
onto $B$. Then by construction for each $k \in \N$ one has that
$F_{B_k}(\Pi(t)_{|B_k})$ is a $\P$ valued EFC process started from
$\1$ and $F_{B_0}(\Pi(t)_{|B_0})$ is a $\P$ valued EFC process
started from $\0.$ Hence, $\Lambda(\pi \cap \Pi(t))$ exists (as
well as the asymptotic frequencies of the blocks of the form
$B_k(t) \cap B_0$).

Now for the second point, for each $j,k \in \N$ define
$$s_{j,k}(n) = \max \{ k \le n : k \in \sigma_{\pi,\pi'}(
B_j(t)\cap B_k)\}.$$

Remark that as $s_{j,k}(n) \nearrow \infty$ when $n \nearrow
\infty$ one has $$\frac{\#\{ k \le
\sigma_{\pi,\pi'}^{-1}(s_{j,k}(n)) : k \in B_j(t) \cap B_k
\}}{\#\{ k \le \sigma_{\pi,\pi'}^{-1}(s_{j,k}(n)) : k \in B_k \}}
\rightarrow \frac{\|B_j(t) \cap B_k\|}{x_{m_{\pi}(k)}}.$$
Furthermore, by definition $$\#\{ k \le n : k \in
\sigma_{\pi,\pi'}(B_k) \} = \#\{ k \le
\sigma_{\pi,\pi'}^{-1}(s_{j,k}(n)) : k \in B_k \}.$$ Hence the
following limit exists and
\begin{eqnarray*}
&& \lim_{n \rightarrow \infty} \frac{1}{n} \#\{ k \le n : k \in
\sigma_{\pi,\pi'}( B_j(t)\cap B_k) \} \\ &=& \lim_{n \rightarrow
\infty} \frac{1}{n} \left(  \frac{\#\{ k \le
\sigma_{\pi,\pi'}^{-1}(s_{j,k}(n)) : k \in B_j(t) \cap B_k
\}}{\#\{ k \le \sigma_{\pi,\pi'}^{-1}(s_{j,k}(n)) : k \in B_k \}}
\right. \\ && \qquad \qquad  \qquad  \#\{ k \le n : k \in
\sigma_{\pi,\pi'}(B_k) \} \Bigg) \\  &=&
x'_{m_{\pi'}(k)}\frac{\|B_j(t) \cap B_k\|}{x_{m_{\pi}(k)}}.
\end{eqnarray*}
The same argument works when $k=0.$ For the last point it is
enough to remark that for each $k$
\begin{eqnarray*}
\|\sigma_{\pi,\pi'}(B_k(t)) \| &=& \| \cup_{i=0}^{\infty}
\sigma_{\pi,\pi'}(B_k(t)\cap B_i)  \| \\ &=& \sum_{i=0}^{\infty}
\|  \sigma_{\pi,\pi'}(B_k(t)\cap B_i)  \|.
\end{eqnarray*}
\end{proof}

The key lemma to prove the proposition is the following:
\begin{lemma}
Consider $\pi_1, \pi_2 \in \P$ with the same hypothesis as in the
above lemma. Suppose furthermore that
$\Lambda(\pi_1)=\Lambda(\pi_2).$ Then
$$Q_{\pi_1} = Q_{\pi_2}.$$
\end{lemma}
\begin{proof}

We have $\pi_1=(B_1^{(1)},B_2^{(1)},...)$ and
$\pi_2=(B_1^{(2)},B_2^{(2)},...).$ Define
$x_1:=\Lambda(\pi_1)=(x_1^{(1)},x_2^{(1)},..)$ and
$x_2:=\Lambda(\pi_2)=(x_1^{(2)},x_2^{(2)},..).$ To ease the
notations, call $\sigma=\sigma_{\pi_1,\pi_2}.$ Note that we have
$\pi_2 =\sigma^{-1}(\pi_1).$

Lemma \ref{exchangeability for evolution} implies that the law of
$\sigma^{-1}(\Pi_1(t))$ is $P_{\pi_2}.$ Lemma \ref{technique 1}
yields that for each $k$ one has
$$\forall t \ge 0 : \Lambda(\sigma(B^{(1)}_k(t)))=\Lambda(B^{(1)}_k(t))$$ and hence
$\Lambda(\sigma^{-1}(\Pi_1(t))=\Lambda(\Pi_1(t)).$ As the
distributions of $\Lambda(\Pi_1(t))$ and
$\Lambda(\sigma^{-1}(\Pi_1(t))$ are respectively $Q_{\pi_1}$ and
$Q_{\pi_2}$ one has $$Q_{\pi_1}=Q_{\pi_2}.$$
\end{proof}

A simple application of Dynkin's criteria (see \cite{dynkin})
concludes the demonstration of the ``Markov" part of Proposition
\ref{Markov pour X}. For the ``Fellerian" part, for $x \in
\sfleche$, call $(Q_x(t),t\ge 0)$ the semi-group of $X$ started
from $X(0)=x.$  As $X$ is right-continuous we must only show that
for $t$ fixed  $x \mapsto Q_x(t)$ is continuous.

Let $x^{(n)} \rightarrow x$ when $n \rightarrow \infty.$ The idea
is to construct a sequence of random variables $X^{(n)}(t)$ each
one with law $Q_{x_n}(t)$ and such that $X^{(n)}(t) \rightarrow
X(t)$ almost surely and where $X$ has law $Q_x(t).$

Take $\pi=(B_0,B_1,B_2,..) \in \P$ such that $\Lambda(\pi)=x.$ For
each $n$ let $\pi_n$ be a partition such that
$\Lambda(\pi_n)=x^{(n)}$ and call $\sigma_n = \sigma_{\pi,\pi_n}$
\footnote{To be rigorous one should extend the definition of
$\sigma_{\pi,\pi_n}$ to allow for the cases where $\pi$ and
$\pi_n$ do not have the same number of blocks.} Furthermore it
should be clear that we can choose $pi_n$ such that for each $k
\le n$ one has $m_{\pi_n}(k) = m_{\pi}(k).$ Hence, one has that
for each $j \ge 0: x^{(n)}_{m_{\pi_n}(j)} \rightarrow
x_{n_{\pi}(j)}$ when $n \rightarrow \infty$ because $x^{(n)}
\rightarrow x.$

As we have remarked, for each $n$ the process $X^{(n)}(t) =
\Lambda((\sigma_n)^{-1}(\Pi(t)))$ where
$\Pi(\cdot)=(B_1(\cdot),B_2(\cdot),..)$ has law $P_{\pi}$ exists
and has law $Q_{x^{(n)}}(t)$.

Using the Lemma \ref{technique 1} one has that
$$\|\sigma_n(B_j(t))\| = \sum_{k \ge 0} \frac{x^{(n)}_{m_{\pi_n}(k)}}{x_{m_{\pi}(k)}} \| B_j(t) \cap B_k \|.$$ This
entails that for each $j$ one has a.s. $$\|\sigma_n(B_j(t))\|
\rightarrow \| B_j(t) \|, \text{ when } n \rightarrow \infty.$$
Hence $Q_{x^{(n)}}(t) \rightarrow Q_x(t)$ in the sense of the
convergence of finite dimensional marginals.

\subsection{Proof of Theorem \ref{propo number}, part 1}
\begin{proof}
We will prove that for each $K \in \N$ one has $\rho(\{\pi : \#\pi
= K\})=0.$

Let us write the equilibrium equations for $\rho^{(n)}(\cdot),$
the invariant measure of the Markov chain $\Pi_{|[n]}$. For each
$\pi \in \P_n$ $$\rho^{(n)} (\pi) \sum_{\pi' \in \P_n \backslash
\{\pi\} } q_n(\pi,\pi') = \sum_{\pi'' \in \P_n \backslash \{ \pi
\}} \rho^{(n)}(\pi'')q_n(\pi'',\pi)$$ where $q_n(\pi,\pi')$ is the
rate at which $\Pi_{|[n]}$ jumps from $\pi$ to $\pi'$. Fix $K \in
\N$ and for each $n \ge K,$ call $A_{n,K} :=\{\pi \in \P_n : \#\pi
\le K\}$ and $D_{n,K}:=\P_n \backslash A_{n,K}$ where $\# \pi$ is
the number of non-empty blocks of $\pi$.

Summing over $A_{n,K}$ yields
\begin{eqnarray*}
&&\sum_{\pi \in A_{n,K}} \rho^{(n)}(\pi)\left[ \sum_{\pi' \in
A_{n,K}\backslash
\{\pi\}}q_n(\pi,\pi') + \sum_{\pi \in D_{n,K}}q_n(\pi,\pi') \right] \\
=&& \sum_{\pi \in A_{n,K}} \left[ \sum_{\pi'' \in
A_{n,K}\backslash \{ \pi \}}\rho^{(n)}(\pi'')q_n(\pi'',\pi) +
\sum_{\pi \in D_{n,K}} \rho^{(n)}(\pi'')q_n(\pi'',\pi)\right]
\end{eqnarray*}
but as
\begin{eqnarray*}
&&\sum_{\pi \in A_{n,K}} \rho^{(n)}(\pi) \left[ \sum_{\pi' \in
A_{n,K}\backslash \{ \pi \}}q_n(\pi,\pi')  \right] \\ = &&
\sum_{\pi \in A_{n,K}} \left[ \sum_{\pi'' \in A_{n,K}\backslash
\{\pi\}}\rho^{(n)}(\pi'')q_n(\pi'',\pi) \right]
\end{eqnarray*}
one has
\begin{eqnarray*}
\sum_{\pi \in A_{n,K}} \rho^{(n)}(\pi)\left[  \sum_{\pi \in
D_{n,K}}q_n(\pi,\pi') \right] = \sum_{\pi \in A_{n,K}} \left[
 \sum_{\pi \in D_{n,K}}
\rho^{(n)}(\pi'')q_n(\pi'',\pi)\right].
\end{eqnarray*}
That is, if we define $q_n(\pi,C)=\sum_{\pi' \in C} q_n(\pi,\pi')$
for each $C \subseteq \P_n,$
\begin{eqnarray} \label{cond d'equilibre}
&&\sum_{\pi \in A_{n,K}} \rho^{(n)}(\pi) q_n(\pi,D_{n,K})  =
\sum_{\pi'' \in D_{n,K}}\rho^{(n)}(\pi'')q_n(\pi'',A_{n,K}).
\end{eqnarray}
Therefore
\begin{eqnarray} \label{cond d'equilibre'}
&&\sum_{\pi \in A_{n,K}\backslash A_{n,K-1}} \rho^{(n)}(\pi)
q_n(\pi,D_{n,K}) \le \sum_{\pi'' \in
D_{n,K}}\rho^{(n)}(\pi'')q_n(\pi'',A_{n,K}). \qquad
\end{eqnarray}

Hence, all we need to prove that $\rho(\{ \pi \in \P :
\#\pi_{|[n]} =K\}) \rightarrow 0$ when $n \rightarrow \infty$ is
to give an upper bound for the right hand-side of (\ref{cond
d'equilibre'}) which is uniform in $n$ and to show that
\begin{eqnarray}\label{converge dutaux min}
\min_{\pi \in A_{n,K} \backslash A_{n,K-1}} q_n(\pi,D_{n,K})
\underset{ n \rightarrow\infty}{\rightarrow} \infty.
\end{eqnarray}

Let us begin with (\ref{converge dutaux min}). Define $$\Phi(q) :=
\ce (q+1) + \int_{\sfleche} (1 - \sum_i x_i^{q+1}) \nuf(dx).$$
This function was introduced by Bertoin in
\cite{bertoin_homogeneous}, where it plays a crucial role as the
Laplace exponent of a subordinator; in particular, $\Phi$ is a
concave increasing function. When $k$ is an integer greater or
equal than 2, $\Phi(k-1)$ is the rate at which $\{[k]\}$ splits,
i.e., it is the arrival rate of atoms $(\pi^{(F)}(t),k(t),t)$ of
$P_F$ such that $\pi^{(F)}_{|[k]}(t) \neq \1_k$ and $k(t)=1.$ More
precisely $\ce k$ is the rate of arrival of atoms that correspond
to erosion and $\int_{\sfleche} (1 - \sum_i x_i^{k}) \nuf(dx)$ is
the rate of arrival of dislocations. Hence, for $\pi \in \P_n$
such that $\#\pi =K$, say $\pi =(B_1,B_2,....,B_K,\o,\o,...)$, one
has
$$q_n(\pi , D_{n,K}) =  \sum_{i : |B_i|>1}
\Phi(|B_i| -1) $$ because it only takes a fragmentation that
creates at least one new block to enter $D_{n,K}.$

First remark that $$\sum_{i:|B_i|>1} \ce |B_i| \ge \ce(n-K+1),$$
next note that $$q \mapsto \int_{\sfleche} (1-\sum_i
x_i^{q+1})\nuf(dx)$$ is also concave and increasing for the same
reason that $\Phi $ is and furthermore
$$ \int_{\sfleche} (1 - \sum_i x_i) \nuf(dx) \ge 0.$$ Hence, for
every $(B_1,...,B_K) \in \P_n$ one has the lower bound
$$ q_n(\pi ,D_{n,K}) = \sum_{i : |B_i|>1} \Phi(|B_i| -1)  \ge
\int_{\sfleche}(1-\sum_i x_i^{(n-K)+1})\nuf(dx) + \ce (n-K+1).$$
As $\Phi(x) \underset{x \rightarrow \infty}{\rightarrow} \infty
\Leftrightarrow \nuf(\sfleche)=\infty \text{ or } \ce >0$ one has
$$\ce >0 \text{ or } \nuf(\sfleche)=\infty \Rightarrow \lim_{n
\rightarrow \infty} \min_{\pi :\# \pi =K} q_n(\pi, D_{n,K}) =
\infty.$$

On the other hand it is clear that $q_n(\pi,A_{n,K})$ only depends
on $\# \pi$ and $K$ (by definition the precise state $\pi$ and $n$
play no role in this rate). By compatibility it is easy to see
that if $\pi,\pi'$ are such that $\# \pi' > \# \pi=K$ then
$$q_n(\pi,A_{n,K}) \ge q_n(\pi',A_{n,K}).$$ Hence, for all $\pi \in D_{n,K}$ one has
$$q_n(\pi,A_{n,K}) \le \tau_K$$ where $\tau_K = q_n(\pi',A_{n,K})$
for all $n$ and any $\pi' \in \P_n$ such that $\# \pi' =K+1,$ and
hence $\tau_K$ is a constant that only depends on $K.$

Therefore
\begin{eqnarray*}
\min_{\pi \in \P_n :\# \pi =K} q_n(\pi, D_{n,K}) \sum_{\pi \in
\P_n :\#\pi =K} \rho^{(n)}(\pi) &\le& \sum_{\pi \in \P_n :\#\pi
=K}
\rho^{(n)}(\pi) q_n(\pi,D_{n,K}) \\
&\le& \sum_{\pi'' \in D_{n,K}}\rho^{(n)}(\pi'')q_n(\pi'',A_{n,K})
\\ &\le&  \tau_K
\sum_{\pi'' \in D_{n,K}}\rho^{(n)}(\pi''),
\end{eqnarray*}
where, on the second inequality, we used (\ref{cond
d'equilibre'}). Thus
$$\rho^{(n)}(\{ \pi \in \P_n : \# \pi =K\}) \le \tau_K / \min_{\pi \in \P_n :\# \pi =K} q_n(\pi,
D_{n,K}).$$ This show that for each $K \in \N,$ one has $\lim_{n
\rightarrow \infty} \rho^{(n)}(\{ \pi \in \P_n : \# \pi =K\})=0$
and thus $\rho(\#\pi<\infty) =0.$
\end{proof}

\subsection{Proof of Theorem \ref{propo number}, part 2}

\begin{proof}
For each $n \in \N$ we define the sequence $(a_i^{(n)})_{i \in
\N}$ by $$a_i^{(n)} := \rho^{(n)}(A_{n,i} \backslash A_{n,i-1})=
\rho^{(n)}(\{ \pi \in \P_n : \# \pi =i \}).$$ We also note $p:=
\nuf(\sfleche)$ the total rate of fragmentation.  The equation
(\ref{cond d'equilibre}) becomes for each $K \in [n]$
\begin{eqnarray} \sum_{\pi : \#\pi =K} \rho^{(n)}(\{\pi\})
q_n(\pi,D_{n,K}) = \sum_{\pi'' \in
D_{n,K}}\rho^{(n)}(\{\pi''\})q_n(\pi'',A_{n,K})
\end{eqnarray}
because the fragmentation is binary. When $\#\pi =K$ one has
$q_n(\pi,D_{n,K}) \le Kp,$ thus
\begin{eqnarray} \label{prop9}
a_K^{(n)} Kp &\ge& \sum_{\pi'' \in
D_{n,K}}\rho^{(n)}(\{\pi''\})q_n(\pi'',A_{n,K})  \notag \\ & \ge &
\sum_{\pi'' : \# \pi'' =K+1}
\rho^{(n)}(\{\pi''\})q_n(\pi'',A_{n,K})
 \notag \\ & \ge & \sum_{\pi'' : \# \pi'' =K+1} \rho^{(n)}(\{\pi''\}) \cc K(K+1)/2
 \notag \\ &\ge & a^{(n)}_{K+1} \cc K(K+1)/2.
\end{eqnarray}
Hence for all $K \in [n-1]$ $$a_K^{(n)} p \ge a_{K+1}^{(n)} \cc
(K+1)/2$$ and thus $$1 = \sum_{i=1}^n a_i^{(n)} <
a_1^{(n)}(1+\sum_{i=1}^{n-1} (p/\cc)^i 2^{i-1}/i !).$$ We conclude
that $a_1^{(n)}$ is uniformly bounded from below by \\
$(1+\sum_{i=1}^{\infty} (p/\cc)^i 2^{i-1}/i!)^{-1}.$ On the other
hand, as $a_1^{(n)}\le 1$ one has
\begin{eqnarray*}
\lim_{n \rightarrow \infty}\sum_{i>K}^n a_i^{(n)} &\le& \lim_{n
\rightarrow \infty} a_1^{(n)} \sum_{i > K-1}^{n-1}
\frac{(2p/\cc)^i}{2i!} \\ &\le & \sum_{i > K-1}^{n-1}
\frac{(2p/\cc)^i}{2i!} \\ &\le & \sum_{i > K-1}^{\infty}
\frac{(2p/\cc)^i}{2i!} \rightarrow 0
\end{eqnarray*}
when $K\rightarrow \infty$. Hence if we define $a_i:= \lim_{n
\rightarrow \infty} a_i^{(n)} = \rho(\{\pi \in \P :\# \pi =i\})$
we have proved that the series $\sum_i a_i$ is convergent and
hence $\lim_{K \rightarrow \infty} \sum_{i
>K} a_i = 0.$ This shows that $\rho(\{ \pi \in
\P :\# \pi =\infty\})=0.$

\end{proof}

\subsection{Proof of Theorem \ref{theorem dust equ}}
\begin{proof}
Define $I_n := \{ \pi=(B_1,B_2,...) \in \P : B_1 \cap [n]=\{1\}
\}$ (when no confusion is possible we sometime use $I_n := \{ \pi
\in \P_n : B_1 =\{1\} \}$) i.e., the partitions of $\N$ such that
the only element of their first block in $[n]$ is $\{1\}$. Our
proof relies on the fact that $$\rho(\{\pi \in \P :
\text{dust}(\pi) \neq \o\}) >0 \Rightarrow \rho(\{\pi \in \P : \pi
\in \cap_n I_n \}) > 0.$$

As above let us write down the equilibrium equations for
$\Pi_{|[n]}(\cdot)$ :
\begin{eqnarray*}
\sum_{\pi \in \P_n \cap I_n} \rho^{(n)}(\pi)q_n(\pi,I_n^c) =
\sum_{\pi' \in I_n^c} \rho^{(n)}(\pi')q_n(\pi,I_n) .
\end{eqnarray*}

Recall that $A_{n,b}$ designates the set of partitions $\pi \in
\P_n$ such that $\# \pi \le b$ and $D_{n,b}=\P_n \backslash
A_{n,b}.$ For each $b$ remark that  $$\min_{\pi \in D_{n,b} \cap
I_n}\{q_n(\pi,I_n^c)\} = q_n(\pi',I_n^c)$$ where $\pi'$ can be any
partition in $\P_n$ such that $\pi' \in I_n$ and $\#\pi'=b+1$. We
can thus define
$$f(b) : = \min_{\pi \in D_{n,b} \cap I_n}\{q_n(\pi,I_n^c)\}.$$
If $\cc>0$ and $\pi \in D_{n,b} \cap I_n$ one can exit from $I_n$
by a coalescence of the Kingman type. This happens with rate
greater than $\cc b.$ If $\nuc(\sfleche)>0$ one can also exit
\textit{via} a coalescence with multiple collision, and this
happens with rate greater than $$ \zeta(b):= \int_{\sfleche}
\left(\sum_i x_i\left(1-\left(1-x_i\right)^{b-1}\right)\right)
\nuc(dx).$$ This $\zeta(b)$ is the rate of arrival of atoms
$\pi^{(C)}(t)$ of $P_C$ such that $\pi^{(C)}(t) \not \in I_b$ and
which do no correspond to a Kingman coalescence. Thus $\sup_{b \in
\N} \zeta(b)$ is the rate of arrival of ``non-Kingman" atoms
$\pi^{(C)}(t)$ of $P_C$ such that $\pi^{(C)}(t) \not \in I:=
\cap_n I_n.$ This rate being $\int_{\sfleche} \left( \sum_i x_i
\right) \nuc(dx)$ and $\zeta(b)$ being an increasing sequence one
has
$$\lim_{b \rightarrow \infty} \zeta(b)= \int_{\sfleche} \left( \sum_i
x_i\right) \nuc(dx).$$ Thus it is clear that, under the conditions
of the proposition, $f(b) \rightarrow \infty$ when $b \rightarrow
\infty.$

On the other hand, when $\pi \in I_n^c$, the rate $q_n(\pi,I_n)$
is the speed at which $1$ is isolated from all the other points,
thus by compatibility it is not hard to see that
$$q_2:= \int_{\sfleche} (1-\sum_i x_i^2)
\nuf(dx) \ge q_n(\pi,I_n) $$ where $q_2$ is the rate at which $1$
is isolated from its first neighbor (the inequality comes from the
inclusion of events).

Hence,
\begin{eqnarray*}
\sum_{\pi \in I_n \cap D_{n,b}}  \rho^{(n)}(\pi) f(b) &\le&
\sum_{\pi \in I_n \cap D_{n,b}}  \rho^{(n)}(\pi)q_n(\pi,I_n^c) \\
&\le& \sum_{\pi' \in I_n^c} \rho^{(n)}(\pi')q_n(\pi,I_n) \\ &\le&
\sum_{\pi' \in I_n^c} \rho^{(n)}(\pi')q_2 \\ &\le & q_2
\end{eqnarray*}
which yields $$ \rho^{(n)}(I_n \cap D_{n,b}) \le q_2 / f(b).$$

Now as $\rho$ is exchangeable one has $\rho(I \cap A_b)=0$ where
$I= \cap_n I_n$ and $A_b=\cap_n A_{n,b}$ (exchangeable partitions
who have dust have an infinite number of singletons, and thus
cannot have a finite number of blocks). Hence $\rho^{(n)}(I_n \cap
A_{n,b}) \rightarrow 0.$

Fix $\epsilon >0$ arbitrarily small and choose $b$ such that $q_2
/ f(b) \le \epsilon /2.$ Then choose $n_0$ such that for all $n
\ge n_0, \rho^{(n)}(I_n \cap A_{n,b}) \le \epsilon /2.$ Hence
$$\forall n \ge n_0 : \rho^{(n)}(I_n) = \rho^{(n)}(I_n \cap
A_{n,b})+\rho^{(n)}(I_n \cap D_{n,b}) \le \epsilon/2 +
\epsilon/2.$$ Thus $\lim_{n \rightarrow \infty} \rho^{(n)}(I_n)=0$
which entails $\rho(B_1=\{1\})=0.$ Now we use the following fact:
$$ \rho(B_1 =\{ 1\}) =\int_{\P} (1-\sum_i \|B_i\|) \rho(d\pi)$$ to
see that $\rho(\text{dust}(\pi)\neq \o)=0.$
\end{proof}

\bibliography{biblio}

\begin{thebibliography}{10}

\bibitem{aldous_survey}
D.~J. Aldous.
\newblock Deterministic and stochastic models for coalescence (aggregation and
  coagulation): a review of the mean-field theory for probabilists.
\newblock {\em Bernoulli}, 5(1):3--48, 1999.

\bibitem{moi_1}
J.~Berestycki.
\newblock Ranked fragmentations.
\newblock {\em ESAIM Probab. Statist.}, 6:157--175 (electronic), 2002.

\bibitem{moi_2}
J.~Berestycki.
\newblock Multifractal spectra of fragmentation processes.
\newblock {\em J. Statist. Phys.}, 113 (3):411--430, 2003.

\bibitem{bertoin_homogeneous}
J.~Bertoin.
\newblock Homogeneous fragmentation processes.
\newblock {\em Probab. Theory Related Fields}, 121(3):301--318, 2001.

\bibitem{bertoin_self-similar}
J.~Bertoin.
\newblock Self-similar fragmentations.
\newblock {\em Ann. Inst. H. Poincar\'e Probab. Statist.}, 38(3):319--340,
  2002.

\bibitem{bertoin_asymptot}
J.~Bertoin.
\newblock The asymptotic behaviour of fragmentation processes.
\newblock {\em J. Euro. Math. Soc.}, 5:395--416, 2003.

\bibitem{cours_bertoin}
J.~Bertoin.
\newblock Fragmentations et coalescences stochastiques.
\newblock In preparation, 2003.

\bibitem{proceeding}
D.~Beysens, X.~Campi, and E.~Peffekorn, editors.
\newblock {\em Proceedings of the workshop : Fragmentation phenomena}, Les
  Houches Series. World Scientific, 1995.

\bibitem{sznitman}
E.~Bolthausen and A.-S. Sznitman.
\newblock On {R}uelle's probability cascades and an abstract cavity method.
\newblock {\em Comm. Math. Phys.}, 197(2):247--276, 1998.

\bibitem{diaconis_mayer-wolf}
P.~Diaconis, E.~Mayer-Wolf, O.~Zeitouni, and M.~P.~W. Zerner.
\newblock Uniqueness of invariant measures for split-merge transformations and
  the poisson-dirichlet law.
\newblock {\em Ann. Probab.}

\bibitem{durrett_limic}
R.~Durrett and V.~Limic.
\newblock A surprising model arising from a species competition model.
\newblock {\em Stoch. Process. Appl}, 102:301--309, 2002.

\bibitem{dynkin}
E.~B. Dynkin.
\newblock {\em Markov processes. {V}ols. {I}}, volume 122 of {\em Translated
  with the authorization and assistance of the author by J. Fabius, V.
  Greenberg, A. Maitra, G. Majone. Die Grundlehren der Mathematischen Wi
  ssenschaften, B\"ande 121}.
\newblock Academic Press Inc., Publishers, New York, 1965.

\bibitem{ethier_kurtz}
S.~N. Ethier and T.~G. Kurtz.
\newblock {\em Markov processes}.
\newblock Wiley Series in Probability and Mathematical Statistics: Probability
  and Mathematical Statistics. John Wiley \& Sons Inc., New York, 1986.
\newblock Characterization and convergence.

\bibitem{jacod_prob_de_mg}
J.~Jacod.
\newblock {\em Calcul stochastique et probl\`emes de martingales}, volume 714
  of {\em Lecture Notes in Mathematics}.
\newblock Springer, Berlin, 1979.

\bibitem{jacod_shiryaev}
J.~Jacod and A.~N. Shiryaev.
\newblock {\em Limit theorems for stochastic processes}, volume 288 of {\em
  Grundlehren der Mathematiscjen Wissenschaften [Fundamental Principles of
  Mathematical Sciences]}.
\newblock Springer-Verlag, Berlin, second edition.

\bibitem{kingman_repre_of_part}
J.~F.~C. Kingman.
\newblock The representation of partition structures.
\newblock {\em J. London Math. Soc. (2)}, 18(2):374--380, 1978.

\bibitem{kingman_genealogy_large_popu}
J.~F.~C. Kingman.
\newblock On the genealogy of large populations.
\newblock {\em J. Appl. Probab.}, (Special Vol. 19A):27--43, 1982.
\newblock Essays in statistical science.

\bibitem{amaury}
A.~Lambert.
\newblock The branching process with logistic growth.
\newblock preprint, 2003.

\bibitem{mohle_sagitov}
M.~M{\"o}hle and S.~Sagitov.
\newblock A classification of coalescent processes for haploid exchangeable
  population models.
\newblock {\em Ann. Probab.}, 29(4):1547--1562, 2001.

\bibitem{pitman_multiple}
J.~Pitman.
\newblock Coalescents with multiple collisions.
\newblock {\em Ann. Probab.}, 27(4):1870--1902, 1999.

\bibitem{pitman_GEM}
J.~Pitman.
\newblock Poisson-{D}irichlet and {GEM} invariant distributions for
  split-and-merge transformation of an interval partition.
\newblock {\em Combin. Probab. Comput.}, 11(5):501--514, 2002.

\bibitem{rogers_william}
L.~C.~G. Rogers and D.~Williams.
\newblock {\em Diffusions, {M}arkov processes, and martingales. {V}ol. 1}.
\newblock Cambridge Mathematical Library. Cambridge University Press,
  Cambridge, 2000.
\newblock Foundations, Reprint of the second (1994) edition.

\bibitem{schweinsberg_multiple}
J.~Schweinsberg.
\newblock Coalescents with simultaneous multiple collisions.
\newblock {\em Electron. J. Probab.}, 5:Paper no.\ 12, 50 pp. (electronic),
  2000.

\bibitem{schweinsberg_comes_down}
J.~Schweinsberg.
\newblock A necessary and sufficient condition for the {$\Lambda$}-coalescent
  to come down from infinity.
\newblock {\em Electron. Comm. Probab.}, 5:1--11 (electronic), 2000.

\end{thebibliography}
\bibliographystyle{abbrv}

\end{document}